\documentclass[twoside,12pt,leqno]{article}
\advance\topmargin by -\headheight\advance\topmargin by -\headsep
\newcommand{\tmvolxx}{xx}
\newcommand{\tmyearyyyy}{yyyy}

\newcommand{\FirstPageHead}[3]{
{\footnotesize 
\vskip -8mm 
\centerline {Travaux math\'ematiques, \quad 
Volume #1 (#2), 
#3,\quad \copyright\  Universit\'e du Luxembourg}}\vspace{-3mm}}

\usepackage{amsmath}
\usepackage{amsthm}
\usepackage{amssymb}
\usepackage{amscd}
\usepackage{graphicx}
\usepackage{epsfig}
\numberwithin{equation}{section}
\newtheorem{theorem}{Theorem}[section]

\theoremstyle{definition}
\newtheorem{definition}[theorem]{Definition}

\newtheorem{remark}[theorem]{Remark}
\newtheorem*{remarku}{Remark}
\numberwithin{equation}{section}

\allowdisplaybreaks
\begin{document}
\thispagestyle{empty}
\FirstPageHead{\tmvolxx}{\tmyearyyyy}{\pageref{firstpage}--\pageref{lastpage}}
\label{firstpage}

\def\KK{K}
\def\cc{c}
\def\zz{z}
\def\empha{\em}
\def\emphas{}
\def\bemphas{}
\def\pemphas{}
\def\scsc{\sc}
\def\sk{k}
\def\scs{\sc}

\newcommand{\acosh}{\mr{acosh}}
\newcommand{\meta}{\zeta}
\newcommand{\halfsum}{\vr}
\newcommand{\weight}{\gamma}
\newcommand{\weightc}{\gamma^c}
\newcommand{\se}{\mr{se}}
\newcommand{\ce}{\mr{ce}}
\newcommand{\vanifac}{f}
\newcommand{\vanisum}{F}
\newcommand{\prin}{1}
\newcommand{\sing}{0}
\newcommand{\pha}{{\mc P}}
\newcommand{\cfg}{{\mc X}}
\newcommand{\BS}{\text{BS}}
\newcommand{\akonst}{c}
\newcommand{\opt}{\mr{opt}}
\newcommand{\Hi}{\mc{H}}
\newcommand{\group}{K}
\newcommand{\lieal}{\mf{k}}
\newcommand{\hol}{\mr{hol}}
\newcommand{\sr}{\mr{S}}
\newcommand{\inco}{\nu}
\newcommand{\tinco}{{\tilde{\nu}}}
\newcommand{\coco}{g}
\newcommand{\vol}{\mr{vol}}
\newcommand{\scapro}[2]{\langle #1,#2 \rangle}
\newcommand{\scaproC}[2]{\langle #1,#2 \rangle}
\newcommand{\scale}{\beta}
\newcommand{\ratio}{\tau}
\newcommand{\vani}{\mc{V}}
\newcommand{\even}{\mr{g}}
\newcommand{\odd}{\mr{u}}
\newcommand{\ket}[1]{{|#1\rangle}}
\newcommand{\abs}[1]{{|#1|}}
\newcommand{\mc}[1]{\mathcal{#1}}
\newcommand{\mr}[1]{\mathrm{#1}}
\newcommand{\SU}{\mr{SU}}
\newcommand{\su}{\mr{su}}
\newcommand{\Gl}{\mr{GL}}
\newcommand{\SO}{\mr{SO}}
\newcommand{\so}{\mr{so}}
\newcommand{\SL}{\mr{SL}}
\newcommand{\ddtn}{\left.\frac{\mr d}{\mr d t}\right|_{t=0}}
\newcommand{\tree}{{\mc T}}
\newcommand{\Pol}{{\mr{Pol}}}
\newcommand{\pr}{{\mr{pr}}}
\renewcommand{\Re}{{\mr{Re}}}
\renewcommand{\Im}{{\mr{Im}}}
\newcommand{\ol}[1]{\overline{#1}}
\newcommand{\tr}{\mathrm{tr}}
\newcommand{\id}{\mathrm{id}}
\newcommand{\RR}{\mathbb{R}}
\newcommand{\CC}{\mathbb{C}}
\newcommand{\ZZ}{\mathbb{Z}}

\newcommand{\II}{\mathbf{1}}

\newcommand{\ctg}{\mr T^\ast}
\newcommand{\tg}{\mr T}
\newcommand{\diag}{{\rm diag}}
\newcommand{\mb}{\mathbf}
\newcommand{\mf}{\mathfrak}
\newcommand{\maxtor}{T}
\newcommand{\U}{{\mr U}}
\newcommand{\Ad}{{\mr{Ad}}}
\newcommand{\stab}{\mr{stab}}
\newcommand{\rref}[1]{{\rm \ref{#1}}}
\newcommand{\vp}{\varphi}
\newcommand{\vr}{\varrho}
\newcommand{\ve}{\varepsilon}
\newcommand{\linie}[3]{\put(#1){\line(#2){#3}}}
\newcommand{\marke}[3]{\put(#1){\put(0.05,0.1){\makebox(-0.1,-0.2)[#2]{$#3$}}}}

\newcommand{\beq}{\begin{equation}}
\newcommand{\eeq}{\end{equation}}
\newcommand{\beqa}{\begin{eqnarray}}
\newcommand{\eeqa}{\end{eqnarray}}
\newcommand{\beqast}{\begin{eqnarray*}}
\newcommand{\eeqast}{\end{eqnarray*}}




\markboth{J.~Huebschmann}{Singular Poisson-K\"ahler geometry and quantization}
$ $
\bigskip

\bigskip

\centerline{{\Large   
Singular Poisson-K\"ahler geometry of stratified K\"ahler 
spaces 
and quantization}}

\bigskip
\bigskip
\centerline{{\large by Johannes Huebschmann}}

\vspace*{.7cm}

\begin{abstract}

In the presence of classical phase space singularities
the standard methods are 
insufficient to attack the problem of quantization.
In certain situations the difficulties 
can be overcome by means of 
K\"ahler quantization on
{\em stratified K\"ahler spaces\/}.
Such a space is a stratified symplectic space together with a complex
analytic structure which is compatible with the stratified symplectic
structure;
in particular each stratum is a K\"ahler manifold in an obvious fashion.
Holomorphic quantization
on a stratified K\"ahler space
then yields a co\-stratified Hilbert space,
a quantum object 
having the classical singularities as its shadow.
Given a K\"ahler manifold with a hamiltonian action of a compact Lie group
that also preserves the complex structure,
reduction after quantization
coincides with quantization after reduction in
the sense that not only the reduced and unreduced quantum phase spaces
correspond but the invariant unreduced and reduced quantum observables as well.

\end{abstract}

\pagestyle{myheadings}

\tableofcontents
\section{Introduction}

In the presence of classical phase space singularities
the standard methods are 
insufficient to attack the problem of quantization.
Ordinary Schr\"odinger quantization leads to a Hilbert space whose 
elements are classes of $L^2$-functions, and incorporating singularities here directly seems 
at present out of sight since we do not know how to handle the singularities in terms of classes of functions. 
However, Hilbert spaces of holomorphic functions are typically spaces
whose points are ordinary functions rather than
classes of functions, and we know well how we can understand
singularities in terms of ordinary functions.
We will show here that, in certain situations, 
by means of a suitable
K\"ahler quantization procedure on
{\em stratified K\"ahler spaces\/}, we can overcome
the difficulties at the quantum level arising from classical phase space singularities.
A stratified K\"ahler space is a stratified symplectic space endowed
with a complex
analytic structure which is compatible with the stratified symplectic
structure;
in particular each stratum is a K\"ahler manifold in an obvious fashion.
K\"ahler quantization then yields a Hilbert space whose points are holomorphic functions
(or more generally holomorphic sections of a holomorphic line bundle);
the resulting quantum Hilbert space actually acquires more structure
which, in turn, has the classical singularities
as its shadow, as we will explain shortly.

Examples of stratified K\"ahler spaces abound:
Symplectic reduction, applied to
K\"ahler manifolds, yields a particular class of examples;
this includes adjoint and generalized adjoint quotients
of complex semisimple Lie groups which,
in turn, underly certain lattice gauge theories.
Other examples come from certain
moduli spaces of holomorphic vector bundles on a Riemann surface
and variants thereof; in physics language, these are
spaces of conformal blocks.
Still other examples arise from 
the closures of holomorphic nilpotent orbits.
Symplectic reduction carries
a  K\"ahler manifold to a  stratified K\"ahler space in such a way
that the sheaf of germs of polarized functions  coincides with the
ordinary sheaf of germs of holomorphic functions.
Projectivization of the closures of
holomorphic nilpotent orbits yields exotic  stratified K\"ahler
structures on complex projective spaces and on certain complex projective
varieties including complex projective quadrics.
Other physical examples are
reduced spaces relative to a constant value of angular momentum.

In the presence of singularities,
a naive approach to quantization might consist
in restriction of the quantization problem 
to a smooth open dense part, the \lq\lq top stratum\rq\rq. 
However this naive procedure
may lead to a loss of information and in fact to
inconsistent results.
To explore the potential impact of classical phase space singularities on quantum problems, 
we developed  the notion of {\em costratified Hilbert space\/}.
This is the appropriate quantum state space over a stratified space; 
a costratified Hilbert space consists of a system of Hilbert spaces, one for each stratum which arises from 
quantization on the closure of that stratum, the stratification is reflected in certain bounded linear 
operators between these Hilbert spaces reversing the partial ordering among the strata, 
and these linear operators are compatible with the quantizations. 
The notion of costratified Hilbert space is, perhaps, the  {\em quantum structure 
having the classical singularities as its shadow\/}.
Within the framework of holomorphic quantization,
a suitable quantization procedure
on stratified K\"ahler spaces
leads to costratified Hilbert spaces.
Given a K\"ahler manifold  with a hamiltonian action of a compact Lie group
that also preserves the complex structure,
reduction after quantization then 
coincides with quantization after reduction in
the sense that not only the reduced and unreduced quantum phase spaces
correspond but the invariant unreduced and reduced quantum observables as well.

We illustrate the approach with a certain concrete model:
In a particular case, we describe a quantum (lattice) gauge theory which
incorporates certain classical singularities. The reduced phase space is
a stratified K\"ahler space, and we make explicit the requisite
singular holomorphic quantization procedure
and spell out the resulting costratified Hilbert space.
In particular, certain tunneling probabilities between
the strata emerge, we will explaine how the energy eigenstates can be determined,
and we will explore corresponding expectation values of the orthoprojectors onto the
subspaces associated with the strata in the strong and  weak
coupling approximations.

The physics described in
the present lecture notes was worked out in research collaboration with 
my physics friends G. Rudoph and M. Schmidt \cite{hurusch}, \cite{varnatwo}. 
I am much indebted to them for having taught me the relevant physics.

\section{Physical systems with classical phase space singularities}

\subsection{An example of a classical phase space singularity}
\label{example}
In $\mathbb R^3$ with coordinates $x,y,r$,
consider the semicone $N$ 
given by the equation $x^2 + y^2 = r^2$ and the inequality $r \geq 0$.
We refer to this semicone as the {\em exotic\/} plane with a single
vertex. The semicone $N$ is the classical reduced
{\em phase space\/} of a single particle 
moving in ordinary affine space of
dimension $\geq 2$ with angular momentum zero.
 In Section \ref{illustration}
below we will  actually justify this claim. The reduced Poisson algebra
$(C^{\infty}N,\{\,\cdot\,,\,\cdot\,\})$ may be described in the following
fashion: Let $x$ and $y$ be the ordinary coordinate functions in
the plane, and consider the algebra $C^{\infty}N$ of smooth
functions in the variables $x,y,r$ subject to the relation $x^2 +
y^2 = r^2$. Define the Poisson bracket $\{\,\cdot\,,\,\cdot\,\}$ on
this algebra by
$$
\{x,y\} = 2r,\  \{x,r\} = 2y,\ \{y,r\} = -2x,
$$
and endow $N$ with the complex structure having $z=x+iy$ as
holomorphic coordinate. The Poisson bracket is then {\em defined
at the vertex\/} as well, 
away from the vertex the Poisson structure is an
ordinary {\em symplectic\/} Poisson structure, and the complex
structure does {\em not\/} \lq\lq see\rq\rq\ the vertex. At the
vertex, the radius function $r$ is {\em not\/} a smooth function
of the variables $x$ and $y$.
Thus the vertex is a singular point for the Poisson structure
whereas it is {\em not\/} a singular point for the complex analytic 
structure.
The Poisson and complex analytic structure
 combine to a \lq\lq stratified K\"ahler structure\rq\rq.
Below we will explain  what this means.

\subsection{Lattice gauge theory}\label{ad}

Let $\group$ be a compact Lie group, let $\lieal$
denote its Lie algebra, and let $\group^{\mathbb C}$ be the
complexification of $\group$. Endow $\lieal$ with an invariant
inner product. The polar decomposition of the complex group
$\group^{\mathbb C}$ and the inner product on $\lieal$ induce a
diffeomorphism
\begin{equation}
\mathrm T^*\group \cong \mathrm T\group \longrightarrow  \group
\times \lieal \longrightarrow \group^{\mathbb C} \label{polar}
\end{equation}
in such a way that the complex structure on $\group^{\mathbb C}$
and the cotangent bundle symplectic structure on $\mathrm
T^*\group$ combine to $\group$-bi-invariant K\"ahler structure.
When we then build a lattice gauge theory from a configuration
space $Q$ which is the product $Q = K^{\ell}$ of $\ell$ copies of
$K$, we arrive at the (unreduced) momentum phase space
\[\mathrm T^*Q= \mathrm T^* K^{\ell} \cong (K^{\mathbb
C})^{\ell},
\]
and reduction modulo the $K$-symmetry given by conjugation leads
to a reduced phase space of the kind
\[
\mathrm T^* K^{\ell}\big/\big/K \cong (K^{\mathbb
C})^{\ell}\big/\big/K^{\mathbb C}
\]
which necessarily involves singularites in a sense to be made
precise, however. 
Here $\mathrm T^* K^{\ell}\big/\big/K$ denotes the symplectic quotient
whereas 
$ (K^{\mathbb
C})^{\ell}\big/\big/K^{\mathbb C}$ refers to the
complex algebraic quotient (geometric invariant theory quotient).
The special case $\ell = 1$, that of a single
spatial plaquette---a quotient of the kind 
$ K^{\mathbb
C}\big/\big/K^{\mathbb C}$ is referred to in the literature
as an {\em adjoint quotient\/}---, 
is mathematically already very attractive and
presents a host of problems which we have elaborated upon in 
\cite{hurusch}. 
To explain how
 in this
particular case
the structure of 
the reduced phase space can be unravelled, 
following \cite{hurusch}, 
we proceed as follows:

Pick a  maximal torus $T$ of $\group$, denote the rank of $T$ by $r$, and
let $W$ be the  Weyl group of  $T$ in $\group$. Then, as a space,
$\mathrm T^*T$ is diffeomorphic to the complexification
$T^{\mathbb C}$ of the torus $T$ and $T^{\mathbb C}$, in turn,
amounts to a product $(\mathbb C^*)^r$ of $r$ copies of the space
$\mathbb C^*$ of non-zero complex numbers. Moreover, the reduced
phase space $\pha$ comes down to the space $\mathrm T^*T\big /W
\cong (\mathbb C^*)^r \big /W$ of $W$-orbits in $(\mathbb C^*)^r$
relative to the action of the Weyl group $W$.

Viewed as the orbit space $\ctg T\big/W$, 
via singular Marsden-Weinstein
reduction,
the reduced
phase space $\pha$ inherits a
stratified symplectic structure. That is to say: (i) The algebra $C^\infty(T^\CC)^W$ of
ordinary smooth $W$-invariant functions on $T^\CC$ inherits a
Poisson bracket and thus furnishes a Poisson algebra of continuous
functions on $\pha$; (ii) for each stratum, the Poisson structure
yields an ordinary symplectic Poisson structure on that stratum;
and (iii) the restriction mapping from $C^{\infty}(T^\CC)^W$ to
the algebra of ordinary smooth functions on that stratum is a
Poisson map.

Viewed as the orbit space $T^\CC\big/W$,  the reduced
phase space $\pha$ acquires a
complex analytic structure. The complex analytic
structure and the Poisson structure combine to a {\em stratified
K\"ahler structure\/} on $\pha$ \cite{kaehler}, \cite{adjoint},
\cite{bedlewo}. The precise meaning of the term \lq\lq
stratified K\"ahler structure\rq\rq\ is that the Poisson structure
satisfies (ii) and (iii) above and that the Poisson and complex
structures satisfy the additional compatibility condition that,
for each stratum, necessarily a complex manifold, the symplectic
and complex structures on that stratum combine to an ordinary
K\"ahler structure.

In Seection \ref{lat} below
we will discuss a model that originates, in
the hamiltonian approach, from lattice gauge theory 
with respect to the group $\group$.
The (classical unreduced) Hamiltonian $H\colon \ctg\group \to \mathbb
R$ of this model is given by
 \begin{equation}
 \label{GHaFn}
H(x,Y)
 =
 -
\frac{1}{2} |Y|^2
 +
\frac{\inco}{2}
 \left(3 - \Re\,\tr(x)\right), \ x\in \group,\, Y \in \lieal\,.
 \end{equation}
Here 
 $\inco = 1/\coco^2$, where $\coco$ is the
coupling constant, the notation
$|\cdot|$ refers to the norm defined by the inner product on
$\lieal$, and the trace refers to some representation of $\group$; below 
we will suppose $\group$ to be realized as a closed subgroup of some unitary group. 
Moreover, the lattice spacing is here set equal
to $1$. The Hamiltonian $H$ is manifestly gauge invariant.

\subsection{The canoe}
\label{cano}
We will now explore the following
special case:
\[
K=\mathrm{SU}(2), \ K^{\mathbb C}=\mathrm{SL}(2,\mathbb C),\
W\cong \mathbb Z/2 .
\]
A maximal torus $T$ in $\mathrm{SU}(2)$ is simply a copy of the
circle group $S^1$, the space $\mathrm T^*T\cong T^{\mathbb C}$ is
a copy of the space $\mathbb C^*$ of non-zero complex numbers, and
the $W$-invariant holomorphic map
\begin{equation}
f \colon \mathbb C^* \longrightarrow \mathbb C,\ f(z) = z+z^{-1}
\label{G-Himap}
\end{equation}
induces a complex analytic isomorphism $\pha\longrightarrow
\mathbb C$ from the reduced space
\[
\pha=\mathrm T^*K \big/\big/K \cong \mathrm T^*T\big /W
\cong\mathbb C^* \big /W
\]
onto a single copy $\mathbb C$ of the complex line. 

\begin{remarku}
More
generally, for $K=\mathrm{SU}(n)$, complex analytically, $\mathrm
T^*K\big/\big/ K$ comes down to $(n-1)$-dimensional complex affine
space  $\mathbb C^{n-1}$. Indeed, $K^{\mathbb
C}=\mathrm{SL}(n,\mathbb C)$, having $(\mathbb C^*)^{n-1}$ as a
maximal complex torus. Realize this torus as the subspace of
$(\mathbb C^*)^n$ which consists of all $(z_1,\dots,z_n)$ such
that $z_1 \dots z_n = 1$. Then the elementary symmetric functions
$\sigma_1, \dots, \sigma_{n-1}$ yield the map
\begin{align*}
(\sigma_1,\dots,\sigma_{n-1})&\colon (\mathbb C^*)^{n-1}
\longrightarrow
\mathbb C^{n-1},\\
\mathbf z=(z_1,\dots,z_n) &\longmapsto (\sigma_1(\mathbf
z),\dots,\sigma_{n-1}(\mathbf z))
\end{align*}
which, in turn, induces the complex analytic isomorphism
\[
\mathrm{SL}(n,\mathbb C) \big/\big/\mathrm{SL}(n,\mathbb C)
\cong(\mathbb C^*)^{n-1} \big /W \cong \mathbb C^{n-1}
\]
from the quotient onto a copy of $\mathbb C^{n-1}$.
We note that, more generally, when $K$ is a general
connected and simply connected Lie group of rank $r$ (say),
in view of an observation of Steinberg's \cite{steinbtw},
the fundamental characters $\chi_1,\ldots,\chi_r$ of $K^{\mathbb C}$
furnish a map from $K^{\mathbb C}$ onto $r$-dimensional
complex affine space $\mathbb A^r$ which identifies the complex
adjoint quotient $K^{\mathbb C}\big/\big/ K^{\mathbb C}$ with 
$\mathbb A^r$.
As a {\em stratified K\"ahler space\/}, the quotient has considerably more structure, though.
We explain this in the sequel for the special case under consideration.
\end{remarku}

Thus we return to
the special case $K=\mathrm{SU}(2)$: In view of the
realization of the complex analytic structure via the holomorphic
map $f \colon \mathbb C^* \to \mathbb C$ given by $f(z) = z+
z^{-1}$ spelled out above, complex analytically, the quotient
$\pha$ is just a copy $\mathbb C$ of the complex line, and we will
take $Z = z+ z^{-1}$ as a holomorphic coordinate on the quotient.
On the other hand, in terms of the notation
\begin{align*}
z&=x+iy, \ Z =X +i Y,\ r^2= x^2+y^2, \\
X &= x + \frac{x}{r^2}, \ Y = y - \frac{y}{r^2},\ \tau =
\frac{y^2}{r^2},
\end{align*}
the real structure admits the following description: In the case
at hand, the algebra written above as $C^\infty(T^\CC)^W$ comes
down the algebra $C^{\infty}(\pha)$ of continuous functions on
$\pha \cong \mathbb C$ which are smooth functions in three
variables (say) $X$, $Y$, $\tau$, subject to certain relations;
the notation $C^{\infty}(\pha)$ is common for such an algebra of
continuous functions even though the elements of this algebra are
not necessarily ordinary smooth functions. To explain the precise
structure of the algebra $C^{\infty}(\pha)$, consider ordinary
real 3-space with coordinates $X$, $Y$, $\tau$ and, in this
3-space, let $C$ be the real semi-algebraic set given by
\[
Y^2 = (X^2 + Y^2 + 4 (\tau -1)) \tau,\quad \tau \geq 0.
\]
As a space, $C$ can be identified with $\pha$. Further, a real
analytic change of coordinates, spelled out in Section 7 of
\cite{bedlewo}, actually identifies $C$ with the familiar {\em
canoe\/}. The algebra $C^{\infty}(\pha)$ is that of Whitney-smooth
functions on $C$, that is, continuous functions on $C$ that are
restrictions of smooth functions in the variables $X$, $Y$, $\tau$
or, equivalently, smooth functions in the variables $X$, $Y$,
$\tau$, where two functions are identified whenever they coincide
on $C$. The Poisson brackets on $C^{\infty}(\pha)$ are determined
by the formulas
\begin{align*}
\{X,Y\} &= X^2 + Y ^2 + 4(2 \tau -1),
\\
\{X,\tau\} &= 2(1- \tau)Y,
\\
\{Y,\tau\} &= 2 \tau X .
\end{align*}
On the subalgebra of $C^{\infty}(\pha)$ which consists of real
polynomial functions in the variables $X$, $Y$, $\tau$, the
relation
\[
Y^2 = (X^2 + Y^2 + 4 (\tau -1)) \tau
\]
is defining. The resulting {\em stratified  K\"ahler\/} structure
on $\pha\cong \mathbb C$ is {\em singular\/} at $-2\in \mathbb C$
and $2\in \mathbb C$, that is, the Poisson structure {\em
vanishes\/} at either of these two points. Further, at $-2\in
\mathbb C$ and $2\in \mathbb C$, the function $\tau$ is {\em
not\/} an ordinary smooth function of the variables $X$ and $Y$,
viz.
\[
\tau = \frac 12 \sqrt{Y^2+ \frac {(X^2 + Y^2-4)^2}{16}} -
\frac{X^2 + Y^2 -4}{8},
\]
whereas away from $-2\in \mathbb C$ and $2\in \mathbb C$, the
Poisson structure is an ordinary symplectic Poisson structure.
This makes explicit, in the case at hand, the singular character
of the reduced space $\pha$ as a stratified K\"ahler space  which,
as a complex analytic space, is just a copy of $\mathbb C$, though
and, as such, has {\em no\/} singularities,  i.~e. is an ordinary
complex manifold.

For later reference, we 
will now describe the stratification of
the reduced configuration space $\cfg \cong T/W$ 
and that of the reduced phase space
$\pha \cong (T\times\mf t)/W$.
The stratifications we will use arise from the $W$-orbit type decompositions:
We will not make precise the notion of stratification and that of stratified space, see e.~g. \cite{gormacon}.

\begin{figure}

 \begin{center}

\unitlength1cm

 \begin{picture}(6,3)
 \put(1,0.27){
 \marke{0,0}{tr}{\pha_+}
 \marke{4,0}{tl}{~\pha_-}
 \marke{4,2}{tl}{~\pha_\prin}
 \put(0,0){\circle*{0.2}}
 \put(4,0){\circle*{0.2}}
 }
 \put(-1,0){
 \put(0,0){\epsfig{file=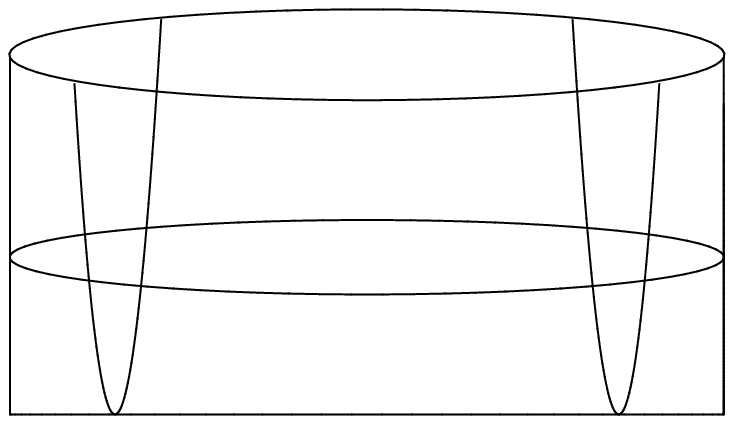,width=6cm,height=3cm}}
  }
 \end{picture}

 \end{center}

 \caption{\label{FKanu} The reduced phase space $\pha$ for
 $\group = \SU(2)$.}

 \end{figure}

The torus $T$ amounts to the complex unit circle
and its Lie algebra $\mf t$ to the imaginary axis. The Weyl group $W=S_2$
acts on $T$ by complex conjugation and on $\mf t$ by reflection.
Hence the reduced configuration space $\cfg \cong T/W$ is homeomorphic
to the closed interval $[-\II,\II]$ and the reduced phase space
$\pha \cong (T\times\mf t)/W$  to the well-known
canoe, see Figure \rref{FKanu}. 

Let 
\[
\cfg_+=\{\II\},\ 
\cfg_-=\{-\II\},\
\cfg_\sing  = \cfg_+\cup\cfg_-=\{-\II,\II\},\ 
\cfg_\prin = ]-\II,\II[
\]
so that the orbit type decomposition of $\cfg$
relative to the $W$-action has the form
$\cfg = \cfg_\prin \cup  \cfg_\sing$.
The \lq\lq piece\rq\rq\ $\cfg_\prin$ (the open interval) is connected;
it is the \lq\lq top\rq\rq\  stratum, the open, connected and dense stratum.
In particular, the restriction to the pre-image of
$\cfg_\prin$
of the orbit projection  is a $W$-covering projection.
The lower stratum $\cfg_\sing$ decomposes into the two connected components
$\cfg_+$ and $\cfg_-$;
the single point in
$\cfg_+$ arises from a fixed point of the $W$-action, and the same is true of $\cfg_-$.
Likewise
the orbit type decomposition of $\pha$
relative to the $W$-action has the form
$\pha = \pha_\prin \cup  \pha_\sing$.
Here  the \lq\lq piece\rq\rq\ 
$\pha_\prin$
 is the \lq\lq top\rq\rq\  stratum, i.~e. the open, connected and dense stratum
which is here 2-dimensional.
As before, the restriction to the pre-image of
$\pha_\prin$ of the orbit projection  is a $W$-covering projection.
Further,
$\pha_\sing$ decomposes into two connected components $\pha_\sing=\pha_+\cup \pha_-$,
each containing a vertex of the canoe;
each such vertex arises from a fixed point of the $W$-action.
Under the identification of $\pha$ with the complex line $\mathbb C$
described previously, the two vertices of the canoe correspond to
the points $2$ and $-2$ of $\mathbb C$ so that
\begin{equation*}
\pha_+ = \{2\}\subseteq \mathbb C,\ \pha_- = \{-2\}\subseteq
\mathbb C, \ \pha_\prin = \mathbb C \setminus \pha_0 =\mathbb C
\setminus \{2,-2\}.
\end{equation*}

A closer look reveals that we can see 
that decomposition  of  $\mathbb C$ as arising from hyperbolic cosine,
viewed as a holomorphic function:
The two points  $2$ and $-2$ are the focal points of the corresponding families of ellipses and hyperbolas in
$\mathbb C$. 
Two of these ellipses and two of these hyperbolas are in fact indicated in Figure 1.
We will come back to the stratifications in Sections \ref{sq} and \ref{lg} below.

\begin{remark}
In the case under discussion ($\group = \SU(2)$), as a stratified symplectic space,
$\pha$ is isomorphic to the reduced phase space of a spherical
pendulum, reduced at vertical angular momentum $0$ (whence the
pendulum is constrained to move in a plane), see \cite{CuBa}.

\end{remark}

\section{Stratified K\"ahler spaces}

In the presence of singularities, restricting quantization to a
smooth open dense stratum, sometimes referred to as
\lq\lq top stratum\rq\rq, can result
in a loss of information and may in fact lead to {\it
inconsistent\/} results.
To develop a satisfactory notion of K\"ahler quantization
in the presence of singularities,
on the classical level, we isolated
a notion of \lq\lq K\"ahler space
with singularities\rq\rq; 
we refer to such a space as
a {\em stratified K\"ahler
space\/}.
Ordinary {\em K\"ahler quantization\/} 
may then be extended to
a {\em quantization scheme over stratified K\"ahler spaces\/}.

We will now explain the concept of a 
\emph{stratified K\"ahler 
space\/}. In \cite{kaehler} we introduced a general
notion of stratified
K\"ahler space 
and that of complex analytic stratified
K\"ahler space as a special case. 
We do not know whether the two notions
really differ. For the present paper, the notion of
complex analytic stratified K\"ahler 
space suffices.
To simplify the terminology somewhat,
\lq\lq stratified K\"ahler space\rq\rq\ will always mean \lq\lq complex
analytic stratified K\"ahler space\rq\rq.

We recall first that, given a stratified space $N$,
a {\em stratified symplectic structure on\/} $N$ is a Poisson algebra
$(C^{\infty}N,\{\,\cdot\,,\,\cdot\,\})$ of continuous
functions on $N$
which, on each stratum, amounts to an ordinary smooth
symplectic Poisson algebra.
The functions in $C^{\infty}N$ are not necessarily
ordinary smooth functions.
Restriction of the functions in $C^{\infty}N$
to a stratum is required to yield
the compactly supported functions on that stratum, and these suffice
to generate a symplectic Poisson algebra on the stratum.

Next we recall that a {\em complex analytic space\/}
(in the sense of {\scsc Grauert\/}) is a topological space $X$, together
with a sheaf of rings $\mathcal O_X$, having the following property:
The space $X$ can be covered by open sets $Y$, each of which
embeds into the polydisc $U$ in some $\mathbb C^n$
(the number $n$ may vary as $U$ varies)
as the zero set of a finite system of holomorphic functions
$f_1,\dots,f_q$ defined on $U$,
such that the 
restriction $\mathcal O_Y$ of the sheaf $\mathcal O_X$ to $Y$
is isomorphic as a sheaf to
the quotient sheaf $\mathcal O_U\big/(f_1,\dots,f_q)$;
here $\mathcal O_U$ is the sheaf of germs of holomorphic functions
on $U$. The sheaf $\mathcal O_X$ is then referred to as the 
{\em sheaf of holomorphic functions on \/} $X$.
See \cite{gunnross} for a development of the general theory of complex 
analytic spaces.

\begin{definition}
 A {\em stratified K\"ahler
space\/} consists of a complex analytic space $N$, together with
\\
{\rm (i)} a complex analytic stratification (a not necessarily proper refinement
of the standard
complex analytic stratification), and with
\\
{\rm (ii)} a stratified symplectic structure
$(C^{\infty}N,\{\,\cdot\,,\,\cdot\,\})$ which is compatible with the complex
analytic structure
\end{definition}
The two structures being {\empha compatible\/}
means the following:
\\
(i) For each point $q$ of $N$ and each holomorphic function
$f$ defined on an open neighborhood $U$ of $q$,
there is an open neighborhood $V$ of $q$ with $V \subset U$
such that, on $V$,
$f$ is the restriction of a function in $C^{\infty}(N)$;
\\
(ii) on each stratum,
the symplectic structure
determined by the symplectic Poisson structure
(on that stratum) combines
with the complex analytic structure to a K\"ahler
structure.

\medskip
\noindent {\scsc Example 1\/}: The {\em exotic plane\/}, endowed with
the structure explained in Subection \ref{example} above, is a
stratified K\"ahler space. Here the radius function $r$ is {\em not\/} an
ordinary smooth function of the variables $x$ and $y$. Thus the stratified
symplectic structure cannot be given in terms of ordinary smooth functions
of the variables $x$ and $y$.

\smallskip
This example generalizes to an entire class of examples:
The {\em closure of a holomorphic nilpotent orbit\/}
(in a hermitian Lie algebra)
inherits a stratified K\"ahler structure
\cite{kaehler}.
Angular momentum zero reduced spaces are special cases thereof;
see Section \ref{illustration} below for details. 

{\em Projectivization\/} of  the closure of a holomorphic
 nilpotent orbit yields what we call an {\em exotic projective variety\/}.
This includes  complex
quadrics, {\scsc Severi\/} and {\scsc Scorza\/} varieties and their
{\em secant\/} varieties \cite{kaehler}, \cite{scorza}.
In physics, spaces of this kind arise as reduced classical phase spaces for
systems of harmonic oscillators with zero angular momentum
and constant energy. We shall explain some of the details in
Section \ref{illustration} below.

\noindent {\scsc Example 2\/}: 
Moduli spaces of semistable holomorphic vector bundles
or, more generally, moduli spaces of semistable principal bundles
on a non-singular complex projective curve
carry  stratified K\"ahler structures
\cite{kaehler}.
These spaces arise
as moduli spaces of homomorphisms
or more generally twisted homomorphisms
from fundamental groups of
surfaces to compact connected Lie groups as well.
In conformal field theory, they occur as spaces of {\em conformal
blocks\/}.
The construction of the moduli spaces
as complex projective varieties goes back to
\cite{narasesh} and \cite{seshaone}; see
\cite{seshaboo} for an exposition of the general theory.
Atiyah and Bott \cite{atibottw}
initiated another approach to the study of these
moduli spaces by identifying them with moduli spaces
of projectively flat constant central curvature
connections on principal bundles over Riemann
surfaces, which they analyzed by methods of gauge theory.
In particular, by applying the method of symplectic
reduction to the action of the infinite-dimensional group
of gauge transformations
on the infinite-dimensional symplectic manifold
of all connections on a principal bundle,
they showed that an invariant inner product on the
Lie algebra of the Lie group in question induces
a natural symplectic structure
on a certain smooth open stratum which, together with
the complex analytic structure,
turns that stratum into an ordinary K\"ahler manifold.
This infinite-dimensional
approach to moduli spaces has  roots in quantum field theory.
Thereafter a finite-dimensional construction of the moduli space
as a symplectic quotient 
arising from an ordinary finite-dimensional
Hamiltonian $G$-space for a compact Lie group $G$
was developed; see  \cite{modus}, \cite{oberwork} and the
literature there; this construction exhibits
the moduli space as a stratified symplectic space.
The stratified K\"ahler structure
mentioned above combines the complex analytic structure
with the stratified symplectic structure;
it includes the K\"ahler manifold structure on the open
and dense stratum.

An important special case is that of
the moduli space of semistable
rank 2 degree zero vector bundles with trivial determinant on a
curve of genus 2. As a space, this is just ordinary
complex projective 3-space, but the 
stratified symplectic structure involves more functions than
just ordinary smooth functions.
The complement of the space of stable vector bundles
is a {\em Kummer surface\/}.
See \cite{locpois}, \cite{oberwork} and the literature there.

Any ordinary K\"ahler manifold is 
plainly a stratified K\"ahler space.
This kind of example generalizes in the following fashion:
For a Lie group $K$, we will
denote its Lie algebra by $\mathfrak k$ and the dual
thereof by 
$\mathfrak k^*$.
The next result says that, roughly speaking,
K\"ahler reduction, applied to an ordinary K\"ahler manifold,
yields  a stratified K\"ahler structure on the reduced space.

\begin{theorem}[\cite{kaehler}] 
\label{kaehler1}
Let $N$ be a K\"ahler manifold, acted upon holomorphically
by a complex Lie group $G$ such that the action, restricted to
a compact real form $K$ of $G$, preserves the
K\"ahler structure and is hamiltonian, with momentum mapping
$\mu \colon N \to \mathfrak k^*$. Then the reduced space
$N_0 = \mu^{-1}(0)\big/K$ inherits a 
stratified K\"ahler structure.
\end{theorem}

For intelligibility, we explain briefly
how the structure on the reduced space
$N_0$ arises.
Details may be found in \cite{kaehler}:
Define $C^{\infty}(N_0)$ to be the 
quotient algebra $C^{\infty}(N)^K\big/I^K$, that is, the algebra
$C^{\infty}(N)^K$ of smooth $K$-invariant functions on $N$, 
modulo the ideal $I^K$ of functions in $C^{\infty}(N)^K$
that vanish on the zero locus $\mu^{-1}(0)$.
The ordinary smooth symplectic Poisson structure $\{\,\cdot\,,\,\cdot\,\}$
on $C^{\infty}(N)$ is $K$-invariant and hence induces a
Poisson structure on the algebra $C^{\infty}(N)^K$ of smooth $K$-invariant
functions on $N$. Furthermore, 
Noether's theorem entails that the ideal $I^K$ is a Poisson ideal,
that is to say, given $f \in 
C^{\infty}(N_0)^K$ and $h\in I^K$, the function 
$\{f,h\}$ is in $I^K$ as well. Consequently the Poisson bracket
$\{\,\cdot\,,\,\cdot\,\}$ descends to a Poisson bracket
$\{\,\cdot\,,\,\cdot\,\}_0$ on
$C^{\infty}(N_0)$.
Relative to the orbit type stratification, the Poisson algebra
$(C^{\infty}N_0,\{\,\cdot\,,\,\cdot\,\}_0)$
turns $N_0$ into a stratified symplectic space.

The inclusion of $\mu^{-1}(0)$ into $N$ passes to a
homeomorphism from
$N_0$ onto the 
categorical $G$-quotient
$N\big/\big/ G$ of $N$ in the category of
complex analytic varieties.
The stratified symplectic structure combines with the
complex analytic structure on
$N\big/\big/ G$ to a stratified K\"ahler structure.
When $N$ is complex algebraic, the complex algebraic $G$-quotient
coincides with the complex analytic $G$-quotient.

Thus, in view of Theorem \ref{kaehler1}, examples of stratified K\"ahler
spaces abound.

\smallskip
\noindent {\scsc Example 3\/}: Adjoint quotients of 
complex reductive Lie groups, see \eqref{ad} above.

\begin{remark}{\rm 
In \cite{atibottw}, {\scsc Atiyah and Bott\/}
raised the issue of \emph{determining the singularities\/}
of moduli spaces of semistable holomorphic vector bundles
or, more generally, of moduli spaces of semistable principal bundles
on a non-singular complex projective curve.
The stratified K\"ahler structure
which we isolated on a moduli space of this kind,
as explained in Example 2 above,
actually determines the singularity structure;
in particular, near any point,  the structure may be understood in terms 
of a suitable local model.
The appropriate notion of singularity
is that of singularity in the sense of stratified
K\"ahler spaces; this notion depends on the entire structure,
not just on the complex analytic structure.
Indeed, the examples spelled out above (the exotic plane
with a single vertex, the exotic plane with two vertices,
the 3-dimensional complex projective space with the Kummer surface
as singular locus, etc.) show
that a point of 
a stratified K\"ahler space may well be a singular point without
being a complex analytic singularity.
}
\end{remark}

\section{Quantum theory and classical singularities}

According to {\scsc Dirac}, the {\em correspondence\/} between a
classical theory and its quantum counterpart should be based on an
analogy between their mathematical structures. 
An interesting
issue is then that of the role of singularities in quantum
problems. Singularities are known to arise in classical phase
spaces. For example, in the hamiltonian picture of a theory,
reduction modulo gauge symmetries leads in general to
singularities on the classical level. This leads to the
question what the significance of singularities on the quantum
side might be. Can we ignore them, or is there a quantum structure
which has the classical singularities as its shadow?
As far as know, one of the first
papers in this topic is that of {\scsc Emmrich and R\"omer\/}
\cite{emmroeme}. This paper indicates that wave functions may
\lq\lq congregate\rq\rq\ near a {\em singular\/} point, which goes
counter to the sometimes quoted statement that {\em singular
points in a quantum problem are a set of measure zero so cannot
possibly be important.\/} In a similar vein, {\scsc Asorey et
al\/} observed that vacuum nodes correspond to the chiral gauge
orbits of reducible gauge fields with non-trivial magnetic
monopole components \cite {asfalolu}. It is also noteworthy that
in classical mechanics and
in classical field theories  singularities in the solution spaces
are the {\em rule rather than the exception\/}. 
This is in particular true for Yang-Mills theories and
for Einstein's gravitational theory where
singularities occur even at some of 
the most interesting and physically relevant solutions, namely at the
symmetric ones.
It is still not
understood what role these singularities might have in quantum
gravity. See, for example, {\scsc Arms, Marsden and Moncrief\/}
\cite{armamonc}, \cite{armamotw} and the literature there.

\section{Correspondence principle and Lie-Rinehart\\ algebras}

To make sense of the {\em correspondence principle\/} in
certain {\em singular\/} situations, one needs a tool which, for the
stratified symplectic
Poisson algebra  on a stratified symplectic space,
serves as a {\em replacement\/} for
the tangent bundle of a smooth symplectic manifold. This
replacement is provided by an appropriate {\em Lie-Rinehart
algebra\/}. This Lie-Rinehart algebra yields in particular a
satisfactory generalization of the Lie algebra of smooth vector fields in the
smooth case. This enables us to put {\em flesh on the bones of
Dirac's correspondence principle in certain singular
situations\/}.

A {\em Lie-Rinehart algebra\/} consists of a commutative
algebra and a Lie algebra with additional structure which
generalizes the mutual structure of interaction between the
algebra of smooth functions and the Lie algebra of smooth vector
fields on a smooth manifold. More precisely:

\begin{definition}
  A {\em Lie-Rinehart\/} algebra
consists of a commutative algebra $A$ and a Lie-algebra $L$ such
that $L$ acts on $A$ by derivations and that $L$ has an $A$-module
structure, and these are required to satisfy
\[
\begin{aligned} {} [\alpha,a\beta]  &= \alpha(a)\beta + a [\alpha,\beta],
\\
(a\alpha)(b) &= a (\alpha(b)),
\end{aligned}
\]
where $a,b \in A$ and $\alpha, \beta \in L$.
\end{definition}

\begin{definition}
An $A$-module $M$ which is also a left $L$-module is called a {\it
left\/} $(A,L)$-module provided
\begin{align} \alpha (a x) &=
\alpha(a) x + a \alpha (x) 
\\
(a \alpha) (x) &=  a (\alpha (x)) 
\end{align}
where  $a \in A, \ x \in M,\ \alpha \in L$.
\end{definition}

We will now explain briefly the Lie-Rinehart algebra associated with a
Poisson algebra; more details may be found in
\cite{poiscoho}, \cite{souriau}, and \cite{lradq}.
Thus, let  $(A,\{\,\cdot\,,\,\cdot\,\})$ be a Poisson algebra.
Let $D_A$ the the $A$-module of formal differentials of $A$ the
elements of which we write as $du$, for $u \in A$.
For $u,v \in A$, the association
$$
(du,dv)\longrightarrow \pi (du,dv) = \{u,v\}
$$
yields an $A$-valued $A$-bilinear skew-symmetric 2-form $\pi=
\pi_{\{\,\cdot\,,\,\cdot\,\}}$ on $D_A$, referred to as
the {\em Poisson\/} 2-{\em form\/} associated with
the Poisson structure $\{\,\cdot\,,\,\cdot\,\}$.
The adjoint
$$
\pi^{\sharp} \colon D_A \longrightarrow \mathrm{Der}(A) =
\mathrm{Hom}_A(D_A,A)
$$
of $\pi$ is a morphism of $A$-modules, and the formula
$$
[a du,b dv] = a\{u,b\} dv + b\{a,v\} du + ab d\{u,v\}
$$
yields a Lie bracket $[\,\cdot\,,\,\cdot\,]$ on $D_A$.

\begin{theorem}[\cite{poiscoho}] The $A$-module structure on $D_A$, the bracket
$[\cdot,\cdot]$, and the morphism
 $\pi^{\sharp}$ of $A$-modules turn the pair $(A,D_A)$
into a Lie-Rinehart algebra.
\end{theorem}

We will write the resulting 
Lie-Rinehart algebra as
$(A,D_{\{\,\cdot\,,\,\cdot\,\}})$.
For intelligibility we recall that,
given a Lie-Rinehart algebra $(A,L)$, 
the Lie algebra $L$ together with the additional 
$A$-module structure 
on $L$ and $L$-module structure on $A$ 
is referred to as an $(\mathbb R,A)$-{\em Lie algebra\/}.
Thus $D_{\{\,\cdot\,,\,\cdot\,\}}$ is an $(\mathbb R,A)$-Lie algebra.

When the Poisson algebra $A$ is the algebra of smooth functions $C^{\infty}(M)$
on a symplectic manifold $M$, 
the $A$-dual $\mathrm{Der}(A) =\mathrm{Hom}_A(D_A,A)$
of $D_A$ amounts to the $A$-module $\mathrm{Vect}(M)$ of smooth vector fields,
and 
\begin{equation}
(\pi^{\sharp},\mathrm{Id})\colon (D_A,A) \longrightarrow
(\mathrm{Vect}(M),C^{\infty}(M))
\label{morphism}
\end{equation}
is a morphism of Lie-Rinehart algebras,
where $(\mathrm{Vect}(M),C^{\infty}(M))$
carries its ordinary Lie-Rinehart structure.
The $A$-module morphism 
$\pi^{\sharp}$ is plainly surjective, and the kernel
consists of those formal differentials which
\lq\lq vanish at each point of\rq\rq\ $M$.

We return to our general Poisson algebra $(A,\{\,\cdot\,,\,\cdot\,\})$.
The Poisson 2-form $\pi_{\{\,\cdot\,,\,\cdot\,\}}$
determines an {\em extension\/}
\begin{equation}
0 
\longrightarrow
A
\longrightarrow 
\overline L_{\{\,\cdot\,,\,\cdot\,\}}
\longrightarrow 
D_{\{\,\cdot\,,\,\cdot\,\}}
\longrightarrow 
0
\label{extension}
\end{equation}
of $(\mathbb R,A)$-Lie algebras which is central 
as an extension of ordinary Lie algebras;
in particular, on the kernel $A$, the Lie bracket is trivial.
Moreover, as $A$-modules,
\begin{equation}
\overline L_{\{\,\cdot\,,\,\cdot\,\}} = A \oplus D_{\{\,\cdot\,,\,\cdot\,\}},
\label{directsum}
\end{equation}
and the Lie bracket on $\overline L_{\{\,\cdot\,,\,\cdot\,\}}$ 
is given by
\begin{equation}
[(a,du),(b,dv)] =
\left(
\{u,b\}+ \{a,v\} - \{u,v\}, d\{u,v\}
\right) ,\quad
a,b,u,v \in A.
\label{liebracket}
\end{equation}
Here we have written \lq\lq $\overline L$\rq\rq\ 
rather than simply $L$ to indicate that
the extension \eqref{extension} represents the {\em negative\/} of the class of
$\pi_{\{\,\cdot\,,\,\cdot\,\}}$
in Poisson cohomology
$\mathrm H_{\mathrm{Poisson}}^2(A,A)$, cf. \cite{poiscoho}.
When $(A,\{\,\cdot\,,\,\cdot\,\})$ is the smooth symplectic Poisson algebra
of an ordinary smooth symplectic manifold,
(perhaps) up to sign, the class of $\pi_{\{\,\cdot\,,\,\cdot\,\}}$
comes essentially 
down to the cohomology class represented by the symplectic structure.

The following concept was introduced in \cite{souriau}.

\begin{definition}
\label{prequantummodule} Given an
$(A\otimes \mathbb C)$-module
$M$,
we refer to
an 
$(A,\overline L_{\{\,\cdot\,,\,\cdot\,\}})$-module
structure 
\begin{equation}
\chi
\colon 
\overline L_{\{\,\cdot\,,\,\cdot\,\}}
\longrightarrow
\mathrm{End}_{\mathbb R}(M)
\label{prequantum}
\end{equation}
on $M$  
as a
{\em prequantum module structure for\/}
$(A,{\{\,\cdot\,,\,\cdot\,\}})$
provided \\
{\rm (i)} the values of $\chi$ lie in
$\mathrm{End}_{\mathbb C}(M)$,
that is to say, 
for $a \in A$ and $\alpha \in 
D_{\{\,\cdot\,,\,\cdot\,\}}$,
the operators $\chi(a,\alpha)$ are complex linear
transformations,
and \\
{\rm (ii)}
for every $a\in A$, with reference to the decomposition \eqref{directsum}, 
we have
\begin{equation}
\chi(a,0) = i\,a\,\mathrm{Id}_M.
\label{complex}
\end{equation}
A pair $(M,\chi)$
consisting of 
an $(A\otimes \mathbb C)$-module $M$ and a prequantum 
module structure
will henceforth be referred to as a {\em prequantum module\/}
(for $(A,\{\,\cdot\,,\,\cdot\,\})$.
\end{definition}

{\em Prequantization\/} now proceeds 
in the following fashion,
cf. \cite{poiscoho}:
The assignment to $a \in A$ of
$(a,da) \in
\overline L_{\{\,\cdot\,,\,\cdot\,\}}$
yields a morphism $\iota$ of real Lie algebras
from
$A$ to
$\overline L_{\{\,\cdot\,,\,\cdot\,\}}$;
thus, for any prequantum module $(M,\chi)$,
the composite of $\iota$ with $-i \chi$
is a representation
$a \mapsto \widehat a$
of the $A$ underlying real Lie algebra
having $M$, viewed as a complex vector space,
as its representation space;
this is a representation by $\mathbb C$-linear operators 
so that any constant acts by multiplication,
that is, 
for any real number $r$,
viewed as a member of $A$,
\begin{equation}
\widehat r = r \,\mathrm{Id}
\label{constants}
\end{equation}
and so that, for $a,b \in A$,
\begin{equation}
\widehat {\{a,b\}} = i\,[\widehat a,\widehat b]
\qquad
\text{(the Dirac condition).}
\label{dirac}
\end{equation}
More explicitly, these operators are given by the formula
\begin{equation}
\widehat a (x) = \frac 1 i \chi(0,da) (x) + ax,
\quad
a \in A,\ x \in M.
\label{prequanti}
\end{equation}
In this fashion, prequantization, that is to say,
the first step in the realization of 
the correspondence principle in one direction,
can be made precise in certain singular situations.

When $(A,\{\,\cdot\,,\,\cdot\,\})$ is the Poisson algebra of
smooth functions on an ordinary smooth sympletic manifold,
this prequantization factors through the morphism
\eqref{morphism} of Lie-Rinehart algebras in such a way that,
on the target, the construction comes down to the ordinary prequantization
construction.

\begin{remarku} In the physics
literature, Lie-Rinehart algebras were explored in a paper by
{\scsc Kastler and Stora\/} under the name {\em Lie-Cartan
pairs\/} \cite{kasstora}.
\end{remarku}

\section{Quantization on stratified K\"ahler spaces}

In the paper \cite{qr} we have shown that
the {\em holomorphic\/} quantization scheme may be extended
to stratified K\"ahler spaces.
We recall the main steps:

\noindent 1) The notion of ordinary K\"ahler polarization
generalizes to that of {\em stratified K\"ahler  polarization\/}.
This concept is defined in terms of the Lie-Rinehart
algebra associated with the stratified symplectic Poisson
structure; it
 specifies {\em polarizations on the
strata\/} and, moreover, encapsulates the {\em mutual positions of
polarizations on the strata\/}.

\noindent Under the circumstances of Theorem \ref{kaehler1},
{\em symplectic reduction
carries a K\"ahler polarization preserved by the symmetries into a
stratified K\"ahler polarization\/}.

\noindent 2) The notion of prequantum bundle generalizes to that
of {\em stratified prequantum module\/}.
Given a stratified K\"ahler space,
a stratified prequantum module is, roughly speaking, a system
of prequantum modules in the sense of 
Definition \ref{prequantummodule}, one for the closure of each stratum,
together with appropriate morphisms among them
which reflect the stratification.

\noindent 3) The notion of quantum Hilbert space
generalizes to that of {\em costratified quantum
Hilbert space\/} in such
a way that the costratified structure reflects the stratification
on the classical level.
{\em Thus the costratified Hilbert space structure is a quantum structure
which has the classical singularities as its shadow.\/}

\noindent 4) The main result says that  $ [Q,R] = 0$, that is,
quantization commutes with reduction \cite{qr}:

\begin{theorem} 
Under the circumstances of {\rm Theorem \ref{kaehler1}},
suppose that the K\"ahler manifold is quantizable
(that is, suppose that the cohomology class of the K\"ahler form
is integral).
When a suitable additional condition
is satisfied, reduction after quantization coincides with
quantization after reduction in the sense that not only the
reduced and unreduced quantum phase spaces correspond but the
(invariant) unreduced and reduced quantum observables as well.
\end{theorem}

What is referred to here as \lq suitable additional condition\rq\ 
is a condition on the behaviour of the gradient flow.
For example, when the K\"ahler manifold is compact, the
condition will automatically be satisfied.

On the reduced level, the resulting classical phase
space involves in general singularities and is a 
stratified K\"ahler space; the appropriate quantum phase space is 
then a costratified Hilbert space.

\section{An illustration arising from angular momentum and holomorphic nilpotent orbits}
\label{illustration}
Let $s$ and $\ell$ be non-zero natural numbers.
The unreduced classical momentum phase space of $\ell$ 
particles in $\mathbb R^s$ 
is real affine space of real dimension $2s \ell$. 
For example, for our solar
system, $s=3$, and $\ell$ is the number of
celestial bodies we take into account, that is, the sun, the planets,
their moons, asteriods, etc., and the true physical phase
space is the reduced space subject to the
(physically reasonable) constraint 
that the total angular momentum of the solar system be constant
and non-zero. 
The shifting trick reduces this case to that of total angular momentum
zero relative to the planar orthogonal group.
The subsequent discussion implies that the reduced phase space relative to the planar
orthogonal group
is the space of complex symmetric $(\ell \times \ell)$-matrices
of rank at most equal to 2.
The true reduced phase space we are looking for then fibers over a semisimple orbit
in $\mathfrak{sp}(\ell,\mathbb R)$
with fiber the space of complex symmetric
$(\ell \times \ell)$-matrices of rank at most equal to 2.
The additional requirement that the total
energy be constant then reduces the system by one more degree of freedom.

We return to the general case. 
Identify real affine space of real dimension $2s \ell$ with the vector space
$(\mathbb R^{2s})^{\times\ell}$ as usual, endow
$\mathbb R^{s}$ with the standard inner product,
$\mathbb R^{2\ell}$ with the standard symplectic structure,
and thereafter
$(\mathbb R^{2s})^{\times\ell}$ with the obvious induced 
inner product and symplectic structure.
The isometry group of the inner product on
$\mathbb R^{s}$ is the orthogonal group
$\mathrm O(s,\mathbb R)$, the group of linear transformations
preserving the symplectic structure on
$\mathbb R^{2\ell}$ is the symplectic group
$\mathrm{Sp}(\ell,\mathbb R)$,
and the actions extend to
 linear $\mathrm O(s,\mathbb R)$- and
$\mathrm{Sp}(\ell,\mathbb R)$-actions
on $(\mathbb R^{2s})^{\times\ell}$ in an obvious
manner. As usual, denote the Lie algebras of
$\mathrm O(s,\mathbb R)$ and
$\mathrm{Sp}(\ell,\mathbb R)$
by
$\mathfrak {so}(s,\mathbb R)$ and
$\mathfrak{sp}(\ell,\mathbb R)$, respectively.

The $\mathrm O(s,\mathbb R)$- and
$\mathrm{Sp}(\ell,\mathbb R)$-actions
on $(\mathbb R^{2s})^{\times\ell}$
are hamiltonian.
To spell out the $\mathrm O(s,\mathbb R)$-momentum mapping
having the value zero at the origin,
identify 
$\mathfrak{so}(s,\mathbb R)$ with its dual
$\mathfrak{so}(s,\mathbb R)^*$ by interpreting
 $a \in\mathfrak{so}(s,\mathbb R)$
as the linear functional on $\mathfrak{so}(s,\mathbb R)$
which assigns $\mathrm{tr}(a {}^t x)$ to 
$x \in\mathfrak{so}(s,\mathbb R)$; here ${}^t x$ refers to the
transpose of the matrix $x$.
We note that, for $s \geq 3$,
\[
(s-2)\mathrm{tr}(a {}^t b) =  -\beta(a,b),\ a,b \in 
\mathfrak{so}(s,\mathbb R),
\]
where $\beta$ is the {\scsc Killing\/} form of
$\mathfrak{so}(s,\mathbb R)$.
Moreover,
for a vector $\mathbf x \in \mathbb R^s$, realized as a column vector,
let ${}^t\mathbf x$ be its transpose, so that ${}^t\mathbf x$ is a row vector.
With these preparations out of the way,
the {\em angular momentum mapping\/} 
\[
 \mu_{\mathrm O} \colon (\mathbb R^{2s})^{\times\ell}
 \longrightarrow \mathfrak{so}(s,\mathbb R)
\]
with reference to the origin is
given by
\[
 \mu_{\mathrm O}(\mathbf q_1,\mathbf p_1,\dots, \mathbf q_\ell,\mathbf p_\ell) 
= \mathbf q_1 {}^t \mathbf p_1 -\mathbf p_1 {}^t \mathbf q_1
+\dots +\mathbf q_\ell {}^t  \mathbf p_\ell 
-\mathbf p_\ell {}^t  \mathbf q_\ell.
\]
Likewise, identify
$\mathfrak{sp}(\ell,\mathbb R)$ 
with its dual $\mathfrak{sp}(\ell,\mathbb R)^*$
 by interpreting
 $a \in\mathfrak{sp}(\ell,\mathbb R)$
as the linear functional on $\mathfrak{sp}(\ell,\mathbb R)$
which assigns $\frac 12 \mathrm{tr}(a x)$ to 
$x \in\mathfrak{sp}(\ell,\mathbb R)$;
we remind the reader that the \emph{Killing\/} form
$\beta$ of $\mathfrak{sp}(\ell,\mathbb R)$
is given by 
\[
\beta(a,b) =2(\ell+1)\mathrm{tr}(ab)
\] 
where
$a,b \in \mathfrak{sp}(\ell,\mathbb R)$.
The $\mathrm{Sp}(\ell,\mathbb R)$-momentum mapping
$$
\mu_{\mathrm {Sp}} \colon (\mathbb R^{2s})^{\times\ell} \longrightarrow
\mathfrak{sp}(\ell,\mathbb R) 
$$
having the value zero at the origin is given by the assignment to
\[
[\mathbf  q_1,\mathbf  p_1, \dots, \mathbf  q_\ell, \mathbf  p_{\ell}]
\in ( \mathbb R^s\times \mathbb R^s)^{\times \ell}
\]
of
\[\left[\begin{array}{cccc}
\left[\mathbf q_j \mathbf p_k\right] & -\left[\mathbf q_j \mathbf q_k\right]\\
\left[\mathbf p_j \mathbf p_k\right] & -\left[\mathbf p_j \mathbf q_k\right]
\end{array}
\right]
\in \mathfrak {sp}(\ell,\mathbb R),
\]
where $\left[\mathbf  q_j \mathbf  p_k\right]$ etc.
denotes the $(\ell \times \ell)$-matrix having the inner products
$\mathbf  q_j \mathbf  p_k$ etc. as entries.

Consider
the $\mathrm O(s,\mathbb R)$-reduced space
$$
N_0 = \mu_{\mathrm O}^{-1}(0)\big/ \mathrm O(s,\mathbb R).
$$
The $\mathrm{Sp}(\ell,\mathbb R)$-momentum mapping
induces an embedding of the reduced space
$N_0$
into $\mathfrak{sp}(\ell,\mathbb R)$.
We now explain briefly  how the image of
$N_0$ in $\mathfrak{sp}(\ell,\mathbb R)$
may be described. More details may be found in \cite{kaehler}, see also 
\cite{scorza}.

Choose a positive complex structure $J$ on 
$\mathbb R^{2\ell}$ which is compatible with $\omega$
in the sense that $\omega(J\mathbf u,J\mathbf v) = 
\omega(\mathbf u,\mathbf v)$ for every 
$\mathbf u, \mathbf v \in \mathbb R^{2\ell}$;
here \lq positive\rq\ means that the associated 
real inner product $\,\cdot\,$ on
$ \mathbb R^{2\ell}$ given by
$\mathbf u \cdot \mathbf v = 
\omega(\mathbf u,J\mathbf v)$ 
for $\mathbf u, \mathbf v \in \mathbb R^{2\ell}$
is positive definite.
The subgroup 
of $\mathrm
{Sp}(\ell,\mathbb R)$
which preserves the complex structure
$J$ is a maximal compact subgroup of $\mathrm
{Sp}(\ell,\mathbb R)$; 
relative to a suitable orthonormal basis,
this group comes down to a copy of the ordinary unitary group
$\mathrm U(\ell)$.
Furthermore, the 
complex structure $J$ induces a {\scsc Cartan\/} decomposition
\begin{equation}
\mathfrak{sp}(\ell,\mathbb R)
 =\mathfrak{u}(\ell) \oplus \mathfrak p;
\label{symp9}
\end{equation}
here 
$\mathfrak
u(\ell)$ is the Lie algebra of $\mathrm U(\ell)$, 
the symmetric constituent $\mathfrak p$ decomposes as the direct sum
\[
\mathfrak p \cong \mathrm S_{\mathbb R}^2[\mathbb
R^\ell]\oplus \mathrm S_{\mathbb R}^2[\mathbb
R^\ell]
\]
of two copies of the real vector space $\mathrm S_{\mathbb R}^2[\mathbb
R^\ell]$ of real symmetric 
$(\ell \times \ell)$-matrices,
and the complex structure $J$ induces 
a complex structure on
$\mathrm S_{\mathbb R}^2[\mathbb
R^\ell]\oplus \mathrm S_{\mathbb R}^2[\mathbb
R^\ell]$ in such a way that
the resulting complex vector space is 
complex linearly
isomorphic
to the complex vector space $\mathrm S_{\mathbb C}^2[\mathbb
C^\ell]$ of complex symmetric 
$(\ell \times \ell)$-matrices
in a canonical fashion.
We refer to a nilpotent orbit $\mathcal O$ in $\mathfrak{sp}(\ell,\mathbb R)$
as being \emph{holomorphic\/}
if the orthogonal projection from $\mathfrak{sp}(\ell,\mathbb R)$
to 
$\mathrm S_{\mathbb C}^2[\mathbb
C^\ell]$, restricted to $\mathcal O$, is a diffeomorphism
from $\mathcal O$ onto its image in $\mathrm S_{\mathbb C}^2[\mathbb
C^\ell]$.
The diffeomorphism from a holomorphic nilpotent orbit
$\mathcal O$ onto its image in $\mathrm S_{\mathbb C}^2[\mathbb
C^\ell]$ 
extends to a homeomorphism from the closure
$\overline {\mathcal O}$ onto its image in
 $\mathrm S_{\mathbb C}^2[\mathbb
C^\ell]$, and the closures of the holomorphic nilpotent orbits
constitute an ascending sequence
\begin{equation}
0 \subseteq \overline {\mathcal O}_1 \subseteq \dots \subseteq
\overline {\mathcal O}_k \subseteq \dots \subseteq \overline
{\mathcal O}_{\ell} \subseteq \mathfrak{sp}(\ell,\mathbb R),\ 
1 \leq k \leq \ell,
\label{symp10}
\end{equation}
such that the orthogonal projection
from $\mathfrak{sp}(\ell,\mathbb R)$ to 
$\mathrm S_{\mathbb C}^2[\mathbb
C^\ell]$, restricted to $\overline{\mathcal O}_{\ell}$, is a homeomorphism
from
$\overline {\mathcal O}_{\ell}$
onto $\mathrm S_{\mathbb C}^2[\mathbb C^\ell]$.
For $1 \leq k \leq \ell$,
this orthogonal projection, 
restricted to $\overline{\mathcal O}_{k}$,
is a homeomorphism from
$\overline{\mathcal O}_{k}$ onto the space
of complex symmetric $(\ell \times \ell)$-matrices of rank at most 
equal to $k$; in particular, each space of the kind
$\overline{\mathcal O}_{k}$ is a \emph{stratified\/} space,
the stratification being given by the rank of the
corresponding complex symmetric 
$(\ell \times \ell)$-matrices.

The Lie bracket of the Lie algebra
$\mathfrak {sp}(\ell,\mathbb R)$
induces a Poisson bracket on the algebra
$C^{\infty}(\mathfrak {sp}(\ell,\mathbb R)^*)$
of smooth functions on the dual
$\mathfrak {sp}(\ell,\mathbb R)^*$
of $\mathfrak {sp}(\ell,\mathbb R)$
in a canonical fashion.
Via the identification
of $\mathfrak {sp}(\ell,\mathbb R)$ with its dual, the Lie bracket
on $\mathfrak {sp}(\ell,\mathbb R)$ induces a Poisson bracket
$\{\,\cdot\,,\,\cdot\,\}$
on  $C^{\infty}(\mathfrak {sp}(\ell,\mathbb R))$. 
Indeed, the assignment to
$a \in \mathfrak {sp}(\ell,\mathbb R)$ of the linear function
\[
f_a \colon \mathfrak {sp}(\ell,\mathbb R) \longrightarrow \mathbb R
\]
given by $f_a(x) = \frac  12 \mathrm{tr}(ax)$
induces a linear isomorphism
\begin{equation}
\mathfrak {sp}(\ell,\mathbb R)
\longrightarrow
\mathfrak {sp}(\ell,\mathbb R)^*;
\label{symp11}
\end{equation}
let
\[
[ \,\cdot\, , \, \cdot \,]^* \colon 
\mathfrak {sp}(\ell,\mathbb R)^*
\otimes \mathfrak {sp}(\ell,\mathbb R)^*
\longrightarrow
\mathfrak {sp}(\ell,\mathbb R)^*
\]
be the bracket on $\mathfrak {sp}(\ell,\mathbb R)^*$
induced by the Lie bracket on $\mathfrak {sp}(\ell,\mathbb R)$.
The Poisson bracket
$\{\,\cdot\,,\,\cdot\,\}$
on the algebra $C^{\infty}(\mathfrak {sp}(\ell,\mathbb R))$
is given by the formula
\[
\{f,h\}(x) = [f'(x), h'(x)]^*(x),\ x \in \mathfrak {sp}(\ell,\mathbb R).
\]
The isomorphism \eqref{symp11}
induces an embedding 
of $\mathfrak {sp}(\ell,\mathbb R)$ into
$C^{\infty}(\mathfrak {sp}(\ell,\mathbb R))$,
and this embedding 
is plainly a morphism 
\[
\delta
\colon 
\mathfrak {sp}(\ell,\mathbb R)
\longrightarrow
C^{\infty}(\mathfrak {sp}(\ell,\mathbb R))
\]
of Lie algebras when 
$C^{\infty}(\mathfrak {sp}(\ell,\mathbb R))$ is viewed
as a real Lie algebra via the Poisson bracket.
In the literature,
a morphism of the kind $\delta$ is referred to as a
\emph{comomentum\/} mapping.

Let $\mathcal O$ be a holomorphic nilpotent orbit.
The embedding of
$\overline{\mathcal O}$
into
$\mathfrak {sp}(\ell,\mathbb R)$ induces
a map from
the algebra $C^{\infty}(\mathfrak {sp}(\ell,\mathbb R))$
of ordinary smooth functions
on $\mathfrak {sp}(\ell,\mathbb R)$
to  the algebra $C^{0}(\overline {\mathcal O})$
of continuous functions on
$\overline {\mathcal O}$,
and we denote the image of
$C^{\infty}(\mathfrak {sp}(\ell,\mathbb R))$
in  $C^{0}(\overline {\mathcal O})$ by
$C^{\infty}(\overline {\mathcal O})$.
By construction, each function in
$C^{\infty}(\overline {\mathcal O})$
is the restriction of an ordinary smooth function on the ambient
space $\mathfrak {sp}(\ell,\mathbb R)$.
Since each stratum of  
$\overline {\mathcal O}$ is an ordinary smooth 
closed submanifold of $\mathfrak {sp}(\ell,\mathbb R)$,
the functions
in $C^{\infty}(\overline {\mathcal O})$, 
restricted to a stratum of 
$\overline {\mathcal O}$,
are ordinary smooth functions on that stratum.
Hence $C^{\infty}(\overline {\mathcal O})$ is a
\emph{smooth structure\/}
on $\overline {\mathcal O}$.
The algebra
$C^{\infty}(\overline {\mathcal O})$
is referred to as the algebra of
{\scsc Whitney\/}-smooth functions on
$\overline {\mathcal O}$,
relative to the embedding of 
$\overline {\mathcal O}$ into the affine space
$\mathfrak {sp}(\ell,\mathbb R)$. 
Under the identification \eqref{symp11},
the orbit $\mathcal O$ passes to a \emph{coadjoint\/}
orbit.
Consequently, under the surjection
$C^{\infty}(\mathfrak {sp}(\ell,\mathbb R)) \to
C^{\infty}(\overline {\mathcal O})$,
the Poisson bracket 
$\{\,\cdot\,,\,\cdot\,\}$
on the algebra $C^{\infty}(\mathfrak {sp}(\ell,\mathbb R))$
descends to a Poisson bracket on
$C^{\infty}(\overline {\mathcal O})$, which we still denote by
$\{\,\cdot\,,\,\cdot\,\}$,
with a slight abuse of notation.
This Poisson algebra turns 
$\overline {\mathcal O}$ into a stratified symplectic space.
Combined with the complex analytic structure coming from the projection
from $\overline {\mathcal O}$ onto the corresponding space
of complex symmetric $(\ell \times \ell)$-matrices,
in this fashion, the space $\overline {\mathcal O}$
acquires a \emph{stratified K\"ahler space\/}
structure.
The composite of the above comomentum mapping $\delta$
with the projection 
from $C^{\infty}(\mathfrak {sp}(\ell,\mathbb R))$ to
$C^{\infty}(\overline {\mathcal O})$
yields
an embedding 
\begin{equation}
\delta_{\mathcal O}
\colon 
\mathfrak {sp}(\ell,\mathbb R)
\longrightarrow
C^{\infty}(\overline {\mathcal O})
\label{delta}
\end{equation}
which is still a morphism of Lie algebras
and therefore a comomentum mapping in the appropriate sense.

The $\mathrm{Sp}(\ell,\mathbb R)$-momentum mapping
induces an embedding of the reduced space
$N_0$
into $\mathfrak{sp}(\ell,\mathbb R)$ which identifies
$N_0$ with the closure
$\overline {\mathcal O}_{\min(s,\ell)}$
of the holomorphic nilpotent orbit
$\mathcal O_{\min(s,\ell)}$ in $\mathfrak{sp}(\ell,\mathbb R)$.
In this fashion, the reduced space
$N_0$ inherits a  stratified K\"ahler structure.
Since the
$\mathrm{Sp}(\ell,\mathbb R)$-momentum mapping
induces an
identification
of $N_0$ with $\overline {\mathcal O}_{s}$
for every $s \leq \ell$ in a compatible manner,
the ascending sequence \eqref{symp10},
and in particular the notion of holomorphic nilpotent orbit,
is actually independent of the choice of complex structure $J$ on 
$\mathrm R^{2 \ell}$. 
For a single particle, i.~e. $\ell = 1$, the description of the reduced space
$N_0$ comes down to that of the semicone given
in Section \ref{example} above.

Thus, when the
number $\ell$ of particles is at most equal to the (real)
dimension $s$ of the space $\mathbb A^s$ in which these particles
move, as a space, the reduced space $N$ amounts to a copy of
complex  affine  space of dimension $\frac {\ell(\ell+1)}2$ and
hence to a copy of real affine space of dimension $\ell(\ell+1)$.
When the
number $\ell$ of particles exceeds the (real)
dimension $s$ of the space in which the particles move,
as a space, the reduced space $N$ amounts to a copy of
the complex  affine  variety of complex symmetric matrices
of rank at most equal to $s$.

\section{Quantization in the situation of the previous class of examples}

In the situation of the previous section,
we will now explain briefly the quantization procedure
developed in \cite{qr}.
Suppose that $s \leq \ell$ (for simplicity), let $m =
s \ell$,
and endow the affine coordinate ring of $\mathbb C^{m}$,
that is, the polynomial algebra $\mathbb C[z_1,\dots,z_m]$,
with the inner product $\,\cdot\,$
given by the standard formula
\begin{equation}
\psi \cdot\psi' = \int \psi \overline{\psi'}
\mathrm e^{-\frac {\mathbf z \overline {\mathbf z}} 2} \varepsilon_m,
\quad
\varepsilon_m= \frac {\omega^m}{(2 \pi)^m m!},
\label{innerprod}
\end{equation}
where $\omega$ refers to the symplectic form on 
$\mathbb C^{m}$.
Furthermore, endow
the polynomial algebra $\mathbb C[z_1,\dots,z_m]$
with the induced $\mathrm O(s,\mathbb R)$-action. 
By construction, the affine complex coordinate ring 
$\mathbb C[\overline {\mathcal O}_s ]$
of $\overline {\mathcal O}_s$
is canonically isomorphic to the algebra
\[
\mathbb C[z_1,\dots,z_m]^{\mathrm O(s,\mathbb R)} 
\]
of $\mathrm O(s,\mathbb R)$-invariants in $\mathbb C[z_1,\dots,z_m]$.
The restriction of the inner product $\,\cdot\,$ 
to $\mathbb C[\overline {\mathcal O}_s ]$
turns 
$\mathbb C[\overline {\mathcal O}_s ]$ into a pre-Hilbert space, and 
{\scsc Hilbert\/}
space completion yields a {\scsc Hilbert\/} space 
which we write as $\widehat{\mathbb C}[\overline {\mathcal O}_s ]$.
This is the Hilbert space which arises by
\emph{holomorphic quantization\/}
on the stratified K\"ahler space
$\overline {\mathcal O}_s$;
see \cite{qr} for details.
On this Hilbert space, the elements of the Lie algebra
$\mathfrak u(\ell)$ of the unitary group $\mathrm U(\ell)$
act in an obvious fashion; indeed,
the elements of 
$\mathfrak u(\ell)$, viewed as functions 
in $C^{\infty}(\overline {\mathcal O}_s)$,
are classical observables which are directly quantizable,
and quantization yields the obvious 
$\mathfrak u(\ell)$-representation on 
$\mathbb C[\overline {\mathcal O}_s ]$.
This construction may be carried out for any $s \leq \ell$
and, for each
$s \leq \ell$,
 the resulting quantizations yields a \emph{costratified
Hilbert space\/} of the kind
\[
\mathbb C  \longleftarrow \widehat{\mathbb C}[\overline {\mathcal O_1}]
\longleftarrow \ldots \longleftarrow \widehat {\mathbb C}[\overline {\mathcal
O_s}].
\]
Here each arrow is just a restriction mapping and is
actually a morphism of representations for the
corresponding quantizable observables, in particular, a morphism
of $\mathfrak u(\ell)$-representations;
each arrow amounts essentially to an orthogonal projection. 
Plainly, the costratified structure
integrates to a costratified $\mathrm U(\ell)$-representation,
i.~e. to a corresponding system of $\mathrm
U(\ell)$-representations. The resulting costratified quantum phase
space for $\overline {\mathcal O_s}$ is a kind of {\em singular\/}
Fock space.
This quantum phase space is
entirely given in terms of {\em data on the reduced level\/}.

Consider the unreduced classical harmonic oscillator
energy $E$ which is given by
$E=z_1 \overline z_1 + \dots +
z_m \overline z_m$; it quantizes to the Euler operator (quantized
harmonic oscillator hamiltonian).
For $s\leq \ell$, 
the reduced classical phase space $Q_s$ of $\ell$ harmonic oscillators in
$\mathbb R^s$ with total angular momentum zero
and fixed energy value which is here encoded in the even number $2k$ fits into an ascending sequence
\begin{equation}
 Q_1 \subseteq \dots \subseteq Q_s \subseteq
\dots \subseteq Q_{\ell} \cong \mathbb C \mathrm P^d
\label{project}
\end{equation}
of stratified K\"ahler spaces where
\[
\mathbb C \mathrm P^d =\mathrm P(\mathrm S^2[\mathbb C^\ell]),
\quad d = \frac {\ell(\ell+1)}2-1.
\]
The sequence \eqref{project} arises from the sequence \eqref{symp10}
by \emph{projectivization\/}.
The parameter $k$ (energy value $2 k$) is
encoded in the Poisson structure. Let $\mathcal O(k)$ be the
$k$'th power of the hyperplane bundle on $\mathbb C \mathrm P^d$,
let
\[
\iota_{Q_s}\colon Q_s \longrightarrow Q_{\ell} \cong \mathbb C
\mathrm P^d
\]
be the inclusion, and let
$
\mathcal O_{Q_s}(k)= \iota_{Q_s}^*\mathcal O(k).
$
The  quantum \emph{Hilbert} space amounts now to the space of holomorphic
sections of $\iota_{Q_s}^*\mathcal O(k)$,
and the resulting {\em costratified quantum Hilbert space\/} has the form
\[
\Gamma^{\mathrm{hol}}(\mathcal O_{Q_1}(k)) \longleftarrow \ldots
\longleftarrow \Gamma^{\mathrm{hol}}(\mathcal O_{Q_s}(k)).
\]
Each vector space $\Gamma^{\mathrm{hol}}(\mathcal O_{Q_{s'}}(k))$
($1 \leq s' \leq s$) is a 
finite-dimensional
representation space for the quantizable
observables in $C^{\infty}(Q_{s})$,
in particular, a $\mathfrak{u}(\ell)$-representation, and this
representation
integrates to a $\mathrm U(\ell)$-representation, and each arrow
is a morphism of representations; similarly as before,
these arrows are just restriction maps.

We will now give a description of the decomposition of the space
\[
\Gamma^{\mathrm{hol}}(\mathcal O_{Q_\ell}(k)) = S_{\mathbb C}^k
[\mathfrak p^*]
\]
of homogeneous degree $k$ polynomial functions on $\mathfrak p
=S_{\mathbb C}^2[\mathbb C^\ell]$ into its \emph{irreducible\/} $\mathrm
U(\ell)$-representations in terms of highest weight vectors.
To this end we note that
coordinates $x_1,\dots,x_\ell$ on $\mathbb C^\ell$ give
rise to coordinates of the kind $\{x_{i,j} = x_{j,i}; \, 1 \leq i,j \leq
\ell\}$ on $S_{\mathbb C}^2 [\mathbb C^{\ell}]$, and the
determinants
\[
 \delta_1 = x_{1,1}, \
\delta_2 = \left | \begin{array}{cccc} x_{1,1}& x_{1,2}\\
                           x_{1,2}& x_{2,2}
                   \end{array}\right|,\ 
\delta_3 = \left | \begin{array}{cccc} x_{1,1}& x_{1,2} & x_{1,3}\\
                           x_{1,2}& x_{2,2} & x_{2,3}\\
                           x_{1,3}& x_{2,3} & x_{3,3}\\
                   \end{array} \right|,
\ \text{etc.}
\]
are highest weight vectors for certain $\mathrm
U(\ell)$-re\-pre\-sen\-ta\-tions. 
For $1 \leq s \leq r$
and $k \geq 1$, the $\mathrm U(\ell)$-representation
$\Gamma^{\mathrm{hol}}(\mathcal O_{Q_s}(k))$ is the sum of the
irreducible representations having as highest weight vectors the
monomials
\[
\delta_1^{\alpha} \delta_2^{\beta} \ldots \delta_s^{\gamma}, \quad
\alpha +2 \beta + \dots + s\gamma = k,
\]
and the restriction morphism
\[
 \Gamma^{\mathrm{hol}}(\mathcal O_{Q_s}(k))\longrightarrow
 \Gamma^{\mathrm{hol}}(\mathcal O_{Q_{s-1}}(k))
\]
has the span of the representations involving $\delta_s$ explicitly
as its kernel and,
restricted to the span of those irreducible representations which 
do \emph{not\/} involve
$\delta_s$, this morphism is an isomorphism.

This situation may be interpreted
in the following fashion: The composite
\[
\mu_{2k}\colon \overline {\mathcal O}_s \subseteq
\mathfrak{sp}(\ell,\mathbb R)\cong\mathfrak{sp}(\ell,\mathbb R)^*
\longrightarrow
\mathfrak u(\ell)^* 
\]
is a singular momentum mapping for the
$\mathrm U(\ell)$-action on
$\overline {\mathcal O}_s$;
actually, the adjoint 
$\mathfrak u(\ell) \to C^{\infty}(\overline {\mathcal O}_s)$
of $\mu^{2k}$
amounts to the composite of \eqref{delta} with the inclusion of
$\mathfrak u(\ell)$ into $\mathfrak {sp}(\ell,\mathbb R)$.
The {\sl irreducible $\mathrm U(\ell)$-representations which
correspond to the coadjoint orbits in the image
\[
\mu_{2k}(O_{s'}\setminus O_{s'-1}) \subseteq \mathfrak u(\ell)^*
\]
of the stratum $O_{s'}\setminus O_{s'-1}$ ($1 \leq s' \leq s$) are
precisely the irreducible representations having as highest weight
vectors the monomials
\[
\delta_1^{\alpha} \delta_2^{\beta} \ldots \delta_{s'}^{\gamma}
\quad (\alpha +2 \beta + \dots + s'\gamma = k)
\]
involving $\delta_{s'}$ explicitly, i.~e. with\/} $\gamma \geq 1$.

\section{Holomorphic half-form quantization on the \\ complexification of
a compact Lie group}

Recall that, given a general compact Lie group $\group$, via the
diffeomorphism \eqref{polar}, the complex structure on
$\group^{\mathbb C}$ and the cotangent bundle symplectic structure
on $\mathrm T^*\group$ combine to $\group$-bi-invariant K\"ahler
structure. A global K\"ahler potential is given by the function
$\kappa$ defined by by
\begin{equation*}
\kappa(x\,\mathrm e^{iY}) =|Y|^2,\ x \in\group,\ Y\in \lieal .
\end{equation*}
The function $\kappa$ being a K\"ahler potential signifies that
the symplectic structure on $\mathrm T^*\group\cong
\group^{\mathbb C}$ is given by $ i
\partial \overline\partial \kappa$. Let $\varepsilon$ denote the
symplectic (or Liouville) volume form on $\mathrm T^*\group\cong
\group^{\mathbb C}$, and let
 $\eta$  be the real
$\group$-bi-invariant (analytic) function on $\group^{\mathbb C}$
given by
\[
\eta(x\,\mathrm e^{iY})
=\sqrt{\left|\frac{\sin(\mathrm{ad}(Y))}{\mathrm{ad}(Y)}\right|},
\ x \in \group, \,Y \in \lieal,
\]
cf. \cite{bhallone} (2.10). Thus $\eta^2$ is the density of Haar
measure on $\group^{\mathbb C}$ relative to Liouville measure
$\varepsilon$.

Half-form K\"ahler quantization on $\group^{\mathbb C}$ leads to
the Hilbert space 
\[
\mathcal HL^2(\group^{\mathbb C},\mathrm
e^{-\kappa/\hbar}\eta \varepsilon)
\] 
of holomorphic functions on
$\group^{\mathbb C}$ that are square-integrable relative to
$\mathrm e^{-\kappa/\hbar}\eta \varepsilon$ \cite{bhallone}. Thus
the scalar product in this Hilbert space is given by
 \begin{equation*}
\scaproC{\psi_1}{\psi_2}
 =
 \frac{1}{\vol(\group)} \int_{\group^\CC}
\ol{\psi_1}\psi_2 \mathrm e^{-\kappa/\hbar}\eta \varepsilon .
 \end{equation*}
Relative to left and right translation, $\mathcal
HL^2(\group^{\mathbb C},\mathrm e^{-\kappa/\hbar}\eta
\varepsilon)$  is a unitary $(\group\times
\group)$-re\-pre\-sen\-ta\-tion, and the Hilbert space associated
with $\pha$ by reduction after quantization is the subspace
\[
\mathcal HL^2(\group^{\mathbb C},\mathrm
e^{-\kappa/\hbar}\eta \varepsilon)^\group
\]
of $\group$-invariants relative to conjugation.

Let $\varepsilon_{T}$ denote the Liouville volume form of $\mathrm
T^*T\cong T^{\mathbb C}$. There is a function $\gamma$ on this
space, made explicit in \cite{hurusch}, such that the restriction
mapping induces an isomorphism
\begin{equation}
\mathcal HL^2(\group^{\mathbb C},\mathrm e^{-\kappa/\hbar}\eta
\varepsilon)^\group
 \longrightarrow
\mathcal H L^2(T^{\mathbb C},\mathrm e^{-\kappa/\hbar} \gamma
\varepsilon_{T})^W
\end{equation}
of Hilbert spaces where the scalar product in $\mc H
L^2(T^{\mathbb C},\mathrm e^{-\kappa/\hbar} \gamma
\varepsilon_{T})^W$ is given by
 \beq\label{GscaproTC}
 \frac{1}{\vol(\group)}
 \int_{T^\CC} \ol{\psi_1}\psi_2 e^{-\kappa/\hbar}\gamma\varepsilon_{T}\,.
 \eeq

\section{Singular quantum structure: costratified \\ Hilbert space}
\label{sq}
Let $N$ be a stratified space. Thus $N$ is a disjoint union
$N=\cup N_{\lambda}$ of locally closed subspaces $N_{\lambda}$,
called {\em strata\/}, each stratum being an ordinary smooth
manifold, and the mutual positions of the strata are made precise
in a way not spelled out here. Let $\mathcal C_N$ be the category
whose objects are the strata of $N$ and whose morphisms are the
inclusions $Y' \subseteq \overline Y$ where $Y$ and $Y'$ range
over strata. We define a {\em costratified Hilbert space\/}
relative to $N$ or associated with the stratification of $N$ to be
a system which assigns a Hilbert space $\mathcal C_Y$ to each
stratum $Y$, together with a bounded linear map $\mathcal C_{Y_2}
\to \mathcal C_{Y_1}$ for each inclusion $Y_1 \subseteq \overline
{Y_2}$ such that, whenever $Y_1 \subseteq \overline {Y_2}$ and
$Y_2 \subseteq \overline {Y_3}$, the composite of $\mathcal
C_{Y_3} \to \mathcal C_{Y_2}$ with $\mathcal C_{Y_2} \to \mathcal
C_{Y_1}$ coincides with the bounded linear map $\mathcal
C_{Y_3}\to\mathcal C_{Y_1}$ associated with the inclusion $Y_1
\subseteq \overline {Y_3}$.

We now explain the construction of the costratified Hilbert space
associated with the reduced phase space $\pha$. This costratified
structure is a  {\em quantum analogue\/} of the {\em orbit type
stratification\/}.

In the Hilbert space
\[\mc H=\mc HL^2(\group^\CC,\mr e^{-\kappa/\hbar}\eta\ve)^\group
\cong \mathcal H L^2(T^{\mathbb C},\mathrm e^{-\kappa/\hbar}
\gamma \varepsilon_{T})^W,
\]
we single out subspaces associated with the strata in an obvious
manner. For the special case
\[K=\mathrm{SU}(2),\  \pha=\mathrm
T^*K\big/\big/ K \cong \mathbb C, \] this comes down to the
following procedure:

The elements of $\mc H$
are ordinary holomorphic functions on $\group^\CC$.
Being $\group$-invariant, they are determined by their
restrictions to $\mathrm T^{\mathbb C}$; these are $W$-invariant
holomorphic functions on $\mathrm T^{\mathbb C}$, and these
$W$-invariant holomorphic functions, in turn, are determined by
the holomorphic functions  on
\[\pha = \group^\CC \big/\big/ \group^\CC \cong T^\CC\big/W\cong
\mathbb C
\]
which they induce on that space.
In terms of the realization of $\pha$ as the complex line $\mathbb C$,
the stratification of $\pha$
reproduced in Subsection \ref{cano} above
 is given by the decomposition $\mathbb C=\pha_+ \cup \pha_- \cup \pha_\prin$ of $\mathbb C$ into
 \begin{equation*}
\pha_+ = \{2\}\subseteq \mathbb C,\ \pha_- = \{-2\}\subseteq
\mathbb C, \ \pha_\prin = \mathbb C \setminus \pha_0 =\mathbb C
\setminus \{2,-2\}.
\end{equation*}
The closed subspaces
\begin{align*}
\vani_{+} &=\{ f \in \mathcal H; f\big|_{\pha_{+}}=0\} \subseteq
\mathcal H
\\
\vani_{-}&=\{ f \in \mathcal H; f\big|_{\pha_{-}}=0\} \subseteq
\mathcal H
\end{align*}
are Hilbert spaces, and we {\em define\/} the Hilbert spaces
$\Hi_{+}$ and $\Hi_{-}$ 
to be the orthogonal complements in
$\mc H$ so that
\[
\mc H =\vani_{+} \oplus \Hi_{+}= \vani_{-} \oplus \Hi_{-};
\]
moreover, we take 
$\Hi_{\prin}$ to be  the entire space $\Hi$.
The resulting system 
\[
\{ \mc H; \ \Hi_{\prin},\ \Hi_{+},\
\Hi_{-}\}, 
\]
together with the corresponding orthogonal projections,
is the {\em costratified Hilbert space\/} associated with the
stratification of $\pha$. By construction, this costratified
Hilbert space structure is a {\em quantum analogue\/} of the {\em
orbit type stratification\/} of $\pha$.

\section{The ho\-lo\-mor\-phic Peter-Weyl theorem}
Choose a dominant Weyl chamber in the maximal torus $\mathfrak t$. Given the highest
weight $\lambda$ (relative to the chosen dominant Weyl chamber),
we will denote by $\chi^{\mathbb C}_{\lambda}$ the irreducible
character of $\group^{\mathbb C}$ associated with $\lambda$.

\begin{theorem}[Holomorphic Peter-Weyl theorem] 
The Hilbert
space 
\[
\mathcal HL^2(K^{\mathbb C},\mathrm e^{-\kappa/\hbar}\eta
\varepsilon)
\] 
contains the vector space $\mathbb C [K^{\mathbb
C}]$ of representative functions on $K^{\mathbb C}$ as a dense
subspace and, as a unitary $(K\times K)$-representation, this
Hilbert space decomposes as the direct sum
$$
 \mathcal HL^2(K^{\mathbb
C},\mathrm e^{-\kappa/\hbar}\eta \varepsilon)\cong \widehat\oplus_{\lambda
\in \widehat{K^{\mathbb C}}} V^*_{\lambda} \otimes V_{\lambda}
$$
of $(K\times K)$-isotypical summands, each such summand being written here as
$V^*_{\lambda} \otimes V_{\lambda}$ where  $V_{\lambda}$ refers to
the irreducible $K$-representation associated with the highest
weight $\lambda$.
\end{theorem}

A proof of this theorem and relevant references can be found in
\cite{holopewe}. The holomorphic Peter-Weyl theorem entails that
the irreducible characters $\chi^{\mathbb C}_{\lambda}$ of
$\group^{\mathbb C}$ constitute a Hilbert space basis of
\[
\mathcal H =\mathcal H L^2(K^{\mathbb C},\mathrm e^{-\kappa/\hbar}
\eta \varepsilon)^K .
\]
Given the highest weight $\lambda$, we will denote by
$\chi_{\lambda}$ the corresponding irreducible character of
$\group$; plainly, $\chi_{\lambda}$ is the restriction to $\group$
of the character $\chi^{\mathbb C}_{\lambda}$. As usual, let
$\rho=\frac 12 \sum_{\alpha\in R^+} \alpha$,  the half sum of the
positive roots and, for a highest weight $\lambda$, let
\begin{equation}
C_{\lambda}:=(\hbar\pi)^{\dim(\group)/2}\mathrm
e^{\hbar|\lambda+\rho|^2},
\label{C}
\end{equation}
where $|\lambda+\rho|$ refers to the norm of $\lambda+\rho$
relative to the  inner product on $\lieal$. In view of the ordinary
Peter-Weyl theorem, the $\{\chi_\lambda\}$'s constitute an {\em
orthonormal} basis of the Hilbert space $L^2(\group,\mr d
x)^\group$.

\begin{theorem}
The assignment to $\chi_{\lambda}$ of
$C_{\lambda}^{-1/2}\chi^{\mathbb C}_{\lambda}$, as $\lambda$
ranges over the highest weights, yields a unitary isomorphism
 \begin{equation}
L^2(\group,\mr d x)^\group
 \longrightarrow
\mc HL^2(\group^\CC,\mr e^{-\kappa/\hbar}\eta\ve)^\group
 \label{H}
 \end{equation}
of Hilbert spaces.
\end{theorem}

By means of this isomorphism, the costratified Hilbert space
structure arising from {\em stratified K\"ahler\/} quantization as
explained earlier carries over to the Schr\"odinger quantization.

\section{Quantum Hamiltonian and Peter-Weyl \\ decomposition} \label{lat}

In the K\"ahler quantization, only the constants are quantizable
while in the \linebreak
Schr\"odinger quantization, functions that are at
most quadratic in generalized momenta are quantizable. In
particular, the classical Hamiltonian 
\eqref{GHaFn}
of our model is quantizable
in the Schr\"odinger quantization, having as associated quantum
Hamiltonian the operator
\begin{equation}
H = -\frac{\hbar^2}{2}\Delta_\group + \frac{\inco}{2} (3-\chi_1)
\label{G-Ham}
\end{equation}
on $L^2(\group,dx)^{\group}$. The operator
$\Delta_\group$, in turn,  arises
from the non-positive Laplace-Beltrami operator
$\tilde\Delta_\group$ associated with the bi-invariant Riemannian
metric on $\group$ as follows: The operator $\tilde\Delta_\group$
is essentially self-adjoint on $C^\infty(\group)$ and has a unique
extension $\Delta_\group$ to an (unbounded) self-adjoint operator
on $L^2(\group,dx)$. The spectrum being discrete, the domain of
this extensions is the space of functions of the form $f=\sum_n
\alpha_n \vp_n$ such that $\sum_n |\alpha_n|^2 \lambda_n^2  <
\infty$ where the $\vp_n$'s range over the eigenfunctions and the
$\lambda_n$'s over the eigenvalues of $\tilde\Delta_\group$.

Since the metric is bi-invariant, so is $\Delta_\group$, whence
$\Delta_\group$ restricts to a self-adjoint operator on
$L^2(\group,\mr d x)^\group$, which we still write as
$\Delta_\group$. By means of the isomorphism \eqref{H}, we then
transfer the Hamiltonian, in particular, the operator
$\Delta_\group$, to a self-adjoint operators on $\mathcal H$.
Schur's lemma then tells us the following:

\noindent (1) Each isotypical $(\group\times \group)$-sum\-mand
$L^2(\group,dx)_\lambda$ of $L^2(\group,dx)$ in the Peter-Weyl
decomposition is an eigenspace  for
$\Delta_\group$;

\noindent (2) the representative functions are eigenfunctions for
$\Delta_\group$;

\noindent (3) the eigenvalue $-\varepsilon_{\lambda}$ of
$\Delta_\group$ corresponding to the highest weight $\lambda$ is
given by
 \begin{equation*}
\varepsilon_{\lambda}=(|\lambda+\rho|^2-|\rho|^2).
 \end{equation*}

Thus, in the holomorphic quantization on $\mathrm T^*K\cong
K^{\mathbb C}$, the free energy operator (i.~e. without potential
energy term) arises as the unique extension of the operator
$-\frac 12 \Delta_K$ on $\mc H$ to an unbounded self-adjoint
operator, and the spectral decomposition thereof refines to the
holomorphic Peter-Weyl decomposition of $\mc H$.

\section{The lattice gauge theory model arising from $\mathrm{SU}(2)$}
\label{lg}

In the rest of the paper we will discuss somewhat informally,
for the special case
where the underlying compact group is $K=\mathrm{SU}(2)$,
some of the implications for the physical interpretation ;
see 
\cite{varnatwo} 
for a leisurely somewhat more complete introduction 
and \cite{hurusch} for a systematic description.

To begin with, we write out the requisite data for the special case
under consideration. We denote the roots of $K=\mathrm{SU}(2)$ 
relative to the dominant Weyl chamber choser earlier
by
$\alpha$ and $-\alpha$, so that $\halfsum = \frac12\alpha$. The invariant inner
product on the Lie algebra $\lieal$ of $K$ is of the form 
 \beq\label{G-scapro}
-\frac{1}{2\scale^2} \tr(Y_1Y_2)
 \,,~~~~~~
Y_1,Y_2\in\lieal\,,
 \eeq
with a scaling factor $\beta>0$ which we will  leave unspecified (e.g.,
$\scale = \frac 1 {\sqrt 8}$ for the Killing form). Then
$$
|\alpha|^2 = 4 \scale^2
 \,,~~~~~~
|\halfsum|^2 = \scale^2
\,.
$$
The highest weights are $\lambda_n = \frac n 2 \alpha$, where $n=0,1,2,\dots$
(twice the spin).  Then
 \beq\label{GEWCn}
\ve_n \equiv \ve_{\lambda_n} = \scale^2 n (n+2)
 \,,~~~~~~
C_n \equiv C_{\lambda_n} = (\hbar\pi)^{3/2} \mr e^{\hbar\scale^2(n+1)^2}\,,
 \eeq
cf. \eqref{C} for the significance of the notation $C_{\lambda_n}$.
We will now write the  complex characters $\chi^\CC_{\lambda_n}$
as $\chi^\CC_n$ ($n \geq 0$).
On $T^\CC$, these complex characters  are given by
 \beq\label{GchiC}
\chi^\CC_n\big(\diag(z,z^{-1})\big)
 =
z^n + z^{n-2} + \cdots + z^{-n}
 \,,~~~~~~
z\in\CC\setminus\{0\}\,,
 \eeq
whereas, on $T$, the corresponding real characters take the form
 \beq\label{Gchi}
\chi_n\big(\diag(\mr e^{\mr i x}, \mr e^{- \mr i x})\big)
 =
\frac{\sin\big((n+1)x\big)}{\sin(x)}
 \,,~~~~~~
x\in \RR\,,~~~ n\geq 0\,.
 \eeq
The Weyl group $W$ permutes the two entries of the elements in $T$. Hence, 
the reduced configuration space $\cfg = T/W$ can be parametrized by $x\in[0,\pi]$ through $x\mapsto 
\diag(\mr e^{\mr i x},\mr e^{-\mr i x})$. In this
parametrization, the measure $v$ on $T$ is given by
$$
v\,\mr d t = \frac{\vol(\group)}{\pi} \sin^2(x) \,\mr d x\,.
$$
It follows that the assignment to $\psi\in C^\infty(T)^W$ of the
function
$$
x\mapsto \sqrt 2 ~\sin x~ \psi\big(\diag(\mr e^{\mr i x},\mr
e^{-\mr i x})\big)
 \,,~~~~~~
x\in[0,\pi]\,,
$$
defines a Hilbert space isomorphism from $L^2(\group,\mr d x)^\group$, realized
as a Hilbert space of $W$-invariant $L^2$-functions on $T$,
onto the ordinary $L^2[0,\pi]$, where the inner product in $L^2[0,\pi]$ is
normalized so that the constant function with value $1$ has norm
$1$. In particular, given $n\geq 0$, the character $\chi_n$ is mapped to the 
function given by the expression
 \beq\label{G-char}
\chi_n(x) = \sqrt 2\, \sin((n+1)x)\,.
 \eeq
In view of the isomorphism between $L^2(\group,\mr d x)^\group$ and $L^2[0,\pi]$
and the isomorphism \eqref{H}, we can work in an
abstract Hilbert space $\mc H$ with a distinguished orthonormal basis
$\{|n\rangle : n=0,1,2,\dots\}$. We achieve the passage to the 
holomorphic realization  
$
\mc H L^2(\group^\CC,\mr e^{-\kappa/\hbar} \eta\ve)^\group,
$
to the 
Schr\"odinger realization $L^2(\group,\mr d x)^\group$, and
to the ordinary $L^2$-realization $L^2[0,\pi]$ by 
substitution of, respectively, 
$C_n^{-1/2}\chi^\CC_n$, $\chi_n$, and $\sqrt 2 \sin(n+1)x$,
for $\ket n$. We remark that plotting wave functions in the realization of
$\mc H$ by $L^2[0,\pi]$ has the advantage that,
directly from 
the graph,
 one can read off  the 
corresponding probability densities with respect to
Lebesgue measure on the parameter space $[0,\pi]$.

We determine the subspaces $\Hi_\tau$ for the special case 
$K=\mathrm{SU}(2)$. The orbit type strata are $\pha_+$, $\pha_-$ and $\pha_1$,
where $\pha_\pm$ consists of the class of $\pm \II$ and $\pha_1 = \pha 
\setminus (\pha_+\cup\pha_-)$. (Recall that via the complex analytic isomorphism 
\eqref{G-Himap}, $\pha_\pm$ is identified with the subset $\{\pm 2\}$ of $\CC$.)
Since $\pha_1$ is dense in $\pha$, the space
$\vani_1$ reduces to zero  and so $\Hi_1 = \Hi$. 
By definition, the subspaces $\vani_+$ and $\vani_-$ consist of the functions 
$\psi\in\Hi$ that satisfy the constraints
 \beq\label{Gvanicond}
\psi(\II) = 0
 \,,~~~~~~
\psi(-\II) = 0
 \,,
 \eeq
respectively. One can check that the system $\{\chi^\CC_n - (n+1) \chi^\CC_0 ~:~
n = 1,2,3,\dots\}$ forms a basis in $\vani_+$ and that the system
$\{\chi^\CC_n + (-1)^n \frac{n+1}{2} \chi^\CC_1 ~:~ n = 0,2,3,\dots\}$
forms a basis in $\vani_-$. Taking the orthogonal complements, we arrive at the
following.

\begin{theorem}\label{PHsing}
The subspaces $\Hi_+$ and $\Hi_-$ have dimension $1$. They are
spanned by the normalized vectors
 \beqa \label{Gpsip}
\psi_+ & := & \frac{1}{N} \sum\nolimits_{n=0}^\infty (n+1)\, \mr
e^{-\hbar\scale^2\,(n+1)^2/2}\, \ket n
 \,,
\\ \label{Gpsim}
\psi_- & := & \frac{1}{N} \sum\nolimits_{n=0}^\infty (-1)^n\, (n+1) \, \mr
e^{-\hbar\scale^2 \, (n+1)^2/2} \, \ket n\,,
 \eeqa
respectively. The normalization factor $N$ is determined by the identity
$$
N^2 = \sum_{n=1}^\infty n^2 \, \mr e^{-\hbar\scale^2 \, n^2}\,.
$$
\end{theorem}

Hence, in Dirac notation, the orthogonal projections
$
\Pi_\pm\colon\Hi\to\Hi_\pm
$
are given by the expressions
 \beq\label{Goproj}
\Pi_\pm
 =
|\psi_\pm\rangle\,\langle\psi_\pm|
 \,.
 \eeq
In terms of the $\theta$-constant
$\theta_3(Q) = \sum_{k=-\infty}^\infty Q^{k^2}$,
the normalization factor $N$ is determined by the identity
 \beq\label{Gnormtheta}
N^2
 =
\frac 1 2 \mr e^{-\hbar\scale^2} \theta_3'(\mr e^{-\hbar\scale^2})\,.
 \eeq
The following figure shows plots of $\psi_\pm$ in the realization of $\Hi$ 
via $L^2[0,\pi]$ for $\hbar\scale^2 = 1/128$ (continuous line), $1/32$ (long
dash), $1/8$ (short dash), $1/2$ (alternating short-long dash).

 \begin{center}

\unitlength1cm

 \begin{picture}(12,2.5)
 \put(0.1,0){
 \marke{2.8,0.5}{cc}{\psi_+}
 \marke{8.6,0.5}{cc}{\psi_-}
 }
 \put(-2.8,-0.2){
 \put(1.8,0){\epsfig{file=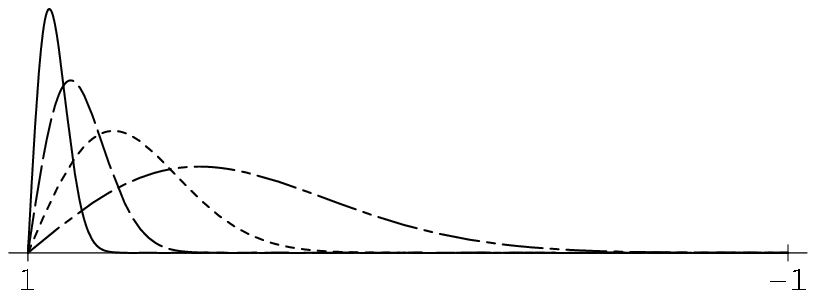,width=6cm,height=2.5cm}}
 \put(7.2,0){\epsfig{file=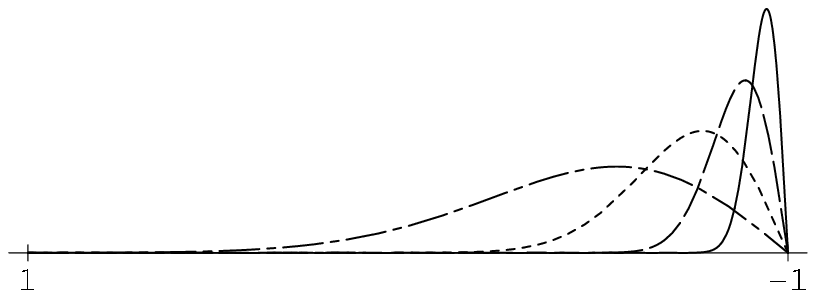,width=6cm,height=2.5cm}}
  }
 \end{picture}

 \end{center}

\section{Tunneling between strata}

Computing the inner product of $\psi_+$ and $\psi_-$,
$$
\langle \psi_+,\psi_-\rangle
 =
 \frac{1}{N^2} \sum_{n=1}^\infty (-1)^{n+1} \,
n^2 \, \mr e^{-\hbar\scale^2\,n^2}
 =
 -\, \frac{
\theta_3'\big(-\mr e^{-\hbar\scale^2}\big)
 }{
\theta_3'\big(\mr e^{-\hbar\scale^2}\big)
 }\,,
$$
we observe that the subspaces $\Hi_+$ and $\Hi_-$ are not orthogonal. They
share a certain overlap which depends on the combined parameter $\hbar\scale^2$.
The absolute square $|\langle \psi_+,\psi_-\rangle|^2$ yields the tunneling
probability between the strata $\pha_+$ and $\pha_-$, i.~e., the probability for
a state prepared at $\pha_+$ to be measured at $\pha_-$ and vice versa. The
following figure shows a plot of the tunneling probability against 
$\hbar \scale^2$. For large values,  this probability tends to $1$ whereas for
$\hbar\scale^2\to 0$, i.e., in the semiclassical limit, it vanishes.
\\

\begin{center}

\hspace*{-2.5cm}\epsfig{file=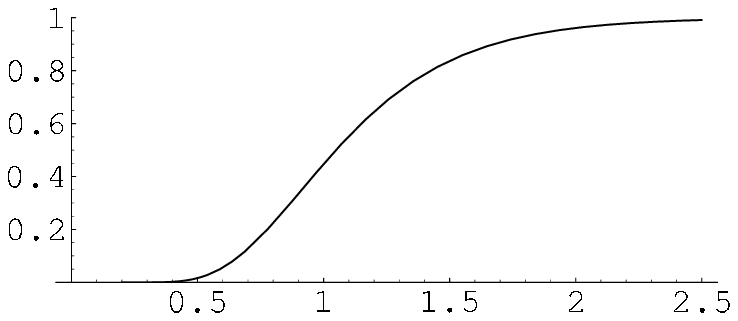,width=8cm}

\end{center}

\section{Energy eigenvalues and eigenstates}

Passing to the realization of $\Hi$ via
$L^2[0,\pi]$ and applying the general formula for the radial part of the
Laplacian on a compact group, see \cite[\S II.3.4]{helgaboo}, 
from the description \eqref{G-Ham} of the quantum Hamiltonian,
viz. 
\[
H = -\frac{\hbar^2}{2}\Delta_\group + \frac{\inco}{2} (3-\chi_1),
\]
we obtain 
the formal expression
$$
 -
\frac{\hbar^2\scale^2}{2}  \left(\frac{\mr d^2}{\mr d x^2} + 1\right)
 +
\frac \nu 2 (3 - \chi_1)
$$
for $H$ on $L^2[0,\pi]$.
Hence
the stationary Schr\"odinger equation can be written
as 
 \beq\label{GSreqf}
 \left(
\frac{\mr d^2}{\mr d x^2}
 +
2 \tinco \cos(x) + \left(\frac{2E}{\hbar^2\scale^2} + 1 -
3\tinco\right)
 \right)
\psi(x)
 =
0\,,
 \eeq
where
 $
\tinco
 =
\frac{\inco}{\hbar^2\scale^2}
 \equiv
\frac{1}{\hbar^2\scale^2\coco^2}
 $,
and $E$ refers to the  eigenvalue. The change of variable $y= (x - \pi)/2$
leads to the Mathieu equation
 \beq\label{GMathieu}
\frac{\mr d^2}{\mr d y^2} f(y) + (a-2q\cos(2y)) f(y) = 0\,,
 \eeq
where
 \beq\label{Gaq}
a = \frac{8 E}{\hbar^2\scale^2} + 4 - 12\tinco
 \,,~~~~~~
q = 4\tinco
 \,;
 \eeq
here $f$ refers to a Whitney smooth function on the interval $[-\pi/2,0]$
satisfying the boundary conditions
 \beq\label{Gboundcond}
f(-\pi/2) = f(0) = 0\,.
 \eeq
For the theory of the Mathieu equation and its solutions, called
{\em Mathieu functions\/}, see
\cite{AbramowitzStegun}. For
certain characteristic values of the parameter $a$ depending
analytically on $q$ and usually denoted by $b_{2n+2}(q)$,
$n=0,1,2,\dots$, solutions satisfying \eqref{Gboundcond} exist.
Given $a = b_{2n+2}(q)$, the corresponding solution is unique up
to a complex factor and can be chosen to be real-valued. It is
usually denoted by $\se_{2n+2}(y;q)$, where \lq$\se$\rq\  stands for
{\em sine elliptic\/}.

Thus, in the realization of $\Hi$ via $L^2[0,\pi]$, the stationary states are
given by
 \beq\label{Gdefxi}
\xi_n(x)
 =
(-1)^{n+1}\sqrt 2
 \left(
\se_{2n+2}\left(\frac{x - \pi}{2};4\tinco\right)
 \right)
 \,,~~~~~~~
n=0,1,2,\dots\,,
 \eeq
and the corresponding eigenvalues by
$$
E_n
 =
\frac{\hbar^2\scale^2}{2}
 \left(
\frac{b_{2n+2}(4\tinco)}{4} + 3 \tinco - 1
 \right)\,.
$$
The factor $(-1)^{n+1}$ ensures that, for $\tinco=0$, we get
$\xi_n = \chi_n$. According to \cite[\S 20.5]{AbramowitzStegun}, for any value
of the parameter $q$, the functions 
$$
\sqrt 2\,\se_{2n+2}(y;q)
 \,,~~~~~~
n = 0,1,2,,\dots\,,
$$
form an orthonormal basis in $L^2[-\pi/2,0]$ and the
characteristic values satisfy $b_2(q)< b_4(q) < b_6(q) < \cdots$. Hence,
the $\xi_n$'s form an orthonormal basis in $\Hi$ and the eigenvalues $E_n$ are
nondegenerate.

\begin{figure} 

\hspace*{-1cm}\epsfig{file=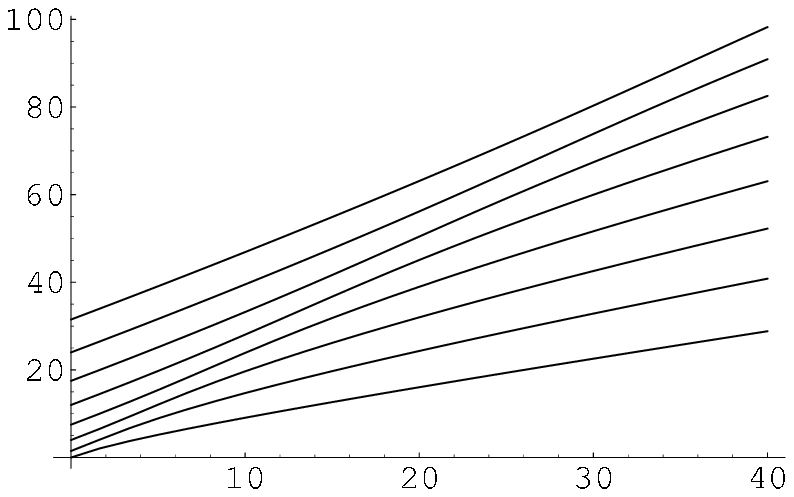,width=6cm}
\hspace*{-0.7cm}\epsfig{file=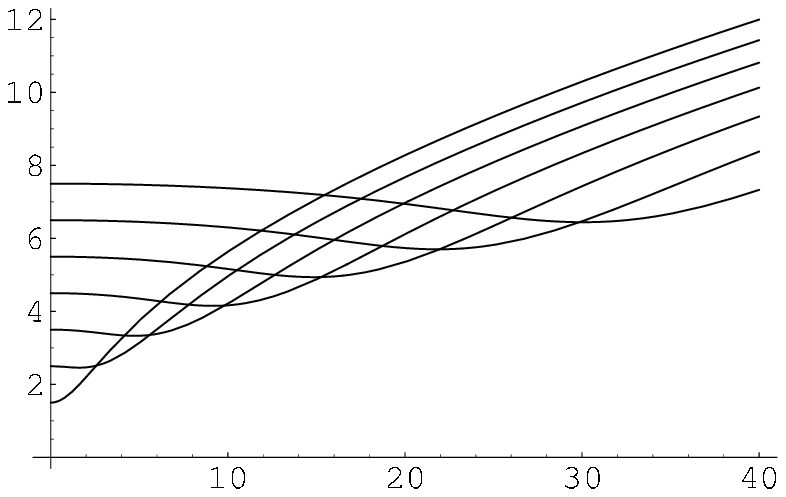,width=6cm}

\caption{\label{AbbEn} Energy eigenvalues $E_n$ (left) and
transition energy values $E_{n+1}-E_n$ (right) for $n=0,\dots,7$ in units of
$\hbar^2\scale^2$ as functions of $\tinco$.}

 \end{figure}

\begin{figure} 

 \begin{center}

 \unitlength1cm

 \begin{picture}(12,5)
 \put(3.05,0){
 \marke{0,2.5}{cc}{\xi_0}
 \marke{5,2.5}{cc}{\xi_1}
 \marke{0,0}{cc}{\xi_2}
 \marke{5,0}{cc}{\xi_3}
 }
 \put(0.25,0){
 \put(0,2.5){\epsfig{file=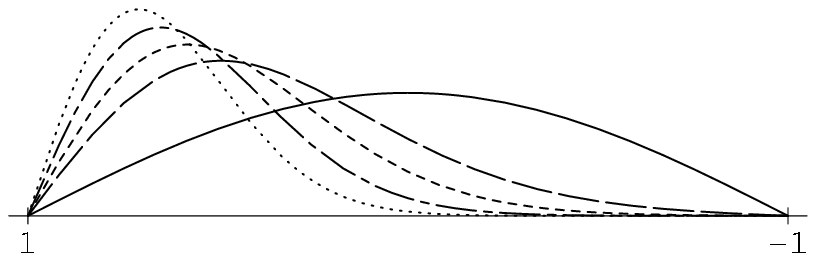,width=4.5cm}}
 \put(5,2.5){\epsfig{file=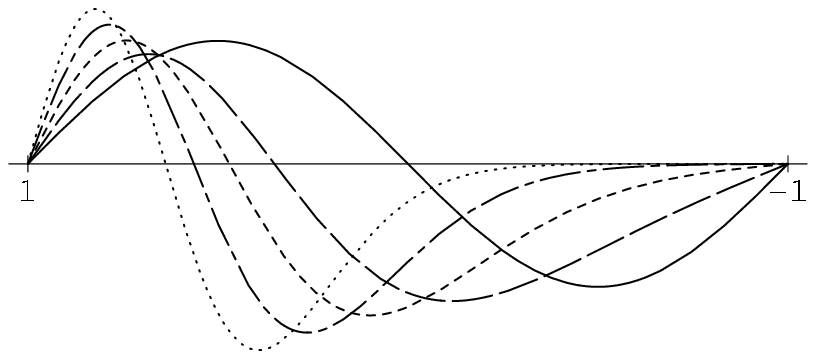,width=4.5cm}}
 \put(0,0){\epsfig{file=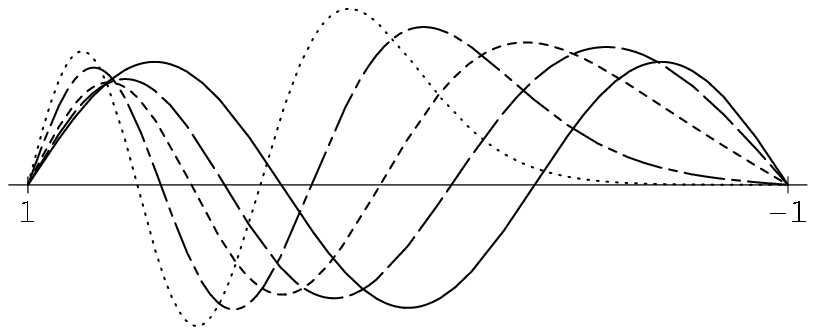,width=4.5cm}}
 \put(5,0){\epsfig{file=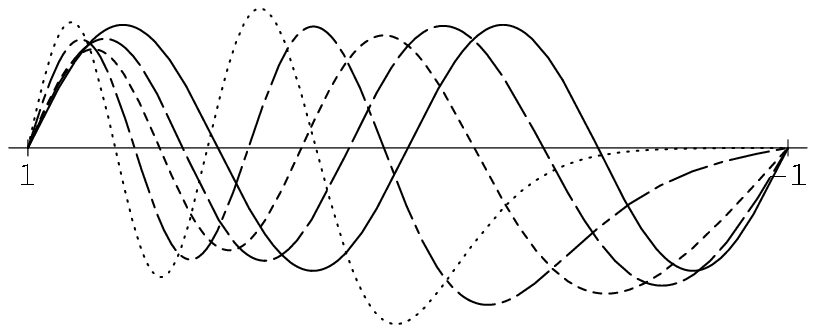,width=4.5cm}}
  }
 \end{picture}

 \end{center}

\caption{\label{Abbxin} Energy eigenfunctions $\xi_0,\dots,\xi_3$ for
$\tinco = 0$ (continuous line), $3$ (long dash), $6$ (short dash),
$12$ (alternating short-long dash), $24$ (dotted line).}

 \end{figure}

Figure \rref{AbbEn} shows the energy eigenvalues $E_n$
and the level separation $E_{n+1}-E_n$ for $n=0,\dots,8$ as
functions of $\tinco$. Figure \rref{Abbxin}
displays the eigenfunctions $\xi_n$, $n=0,\dots,3$, for $\tinco = 0,3,6,12,24$.
The plots have been generated by means of the built-in Mathematica functions
{\tt MathieuS} and {\tt MathieuCharacteristicB}.

\section{Expectation values of the costratification orthoprojectors}

\label{Sexpect}

\begin{figure}[h]

\unitlength1cm

 \begin{picture}(10,11.5)
 \put(-1.5,9){
 \put(0.5,0){\epsfig{file=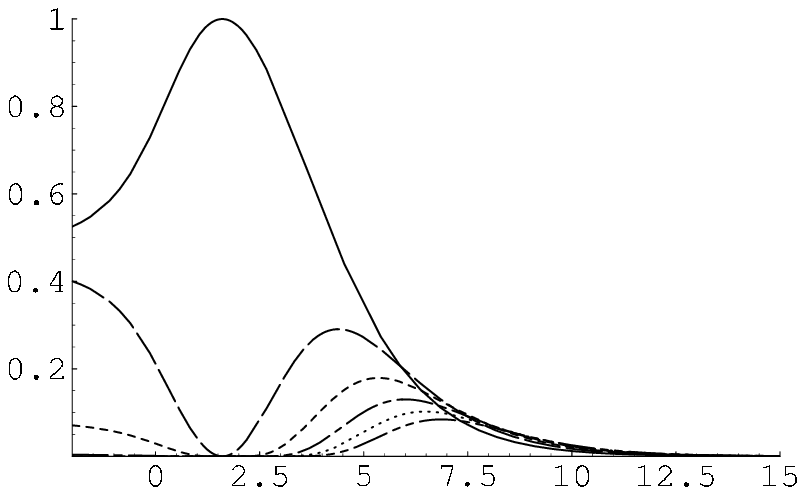,width=6cm,height=2.5cm}}
 \put(6,0){\epsfig{file=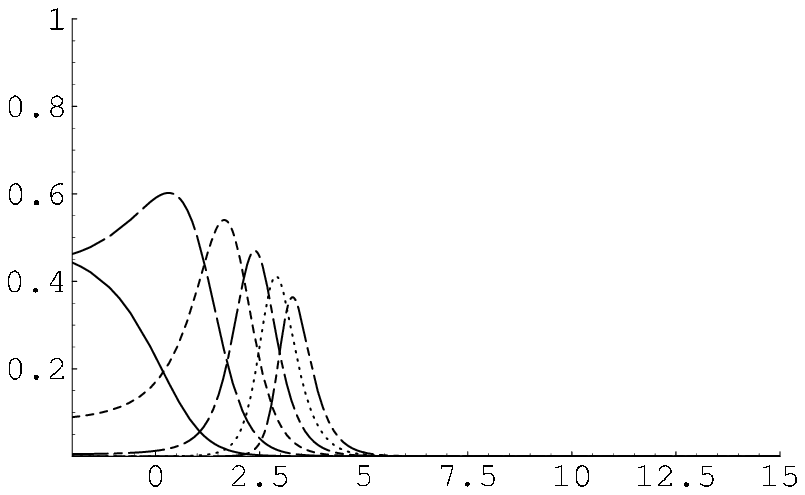,width=6cm,height=2.5cm}}
  }
 \put(-1.5,5){
 \put(0.5,0){\epsfig{file=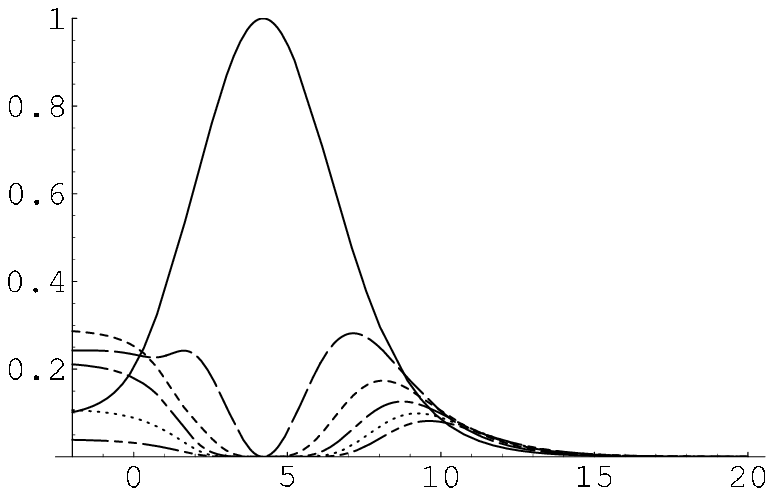,width=6cm,height=2.5cm}}
 \put(6,0){\epsfig{file=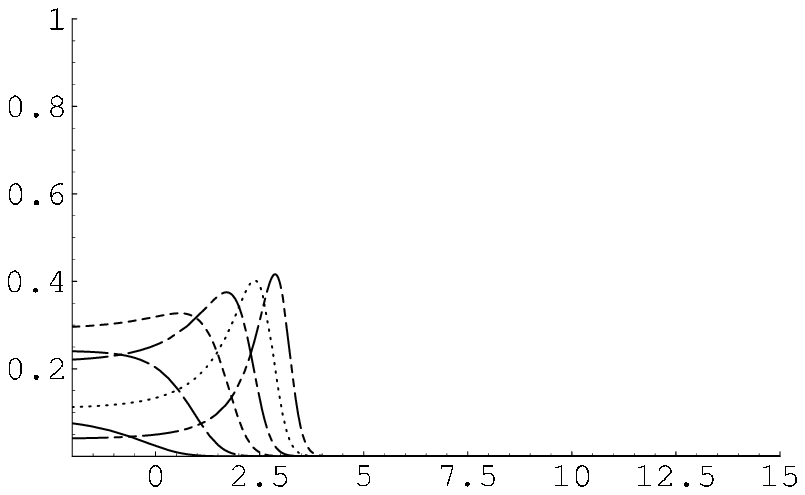,width=6cm,height=2.5cm}}
  }
 \put(-1.5,1){
 \put(0.5,0){\epsfig{file=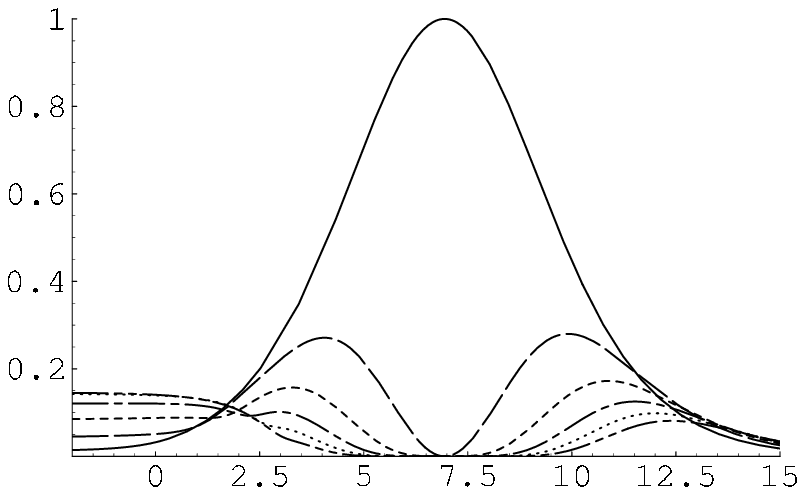,width=6cm,height=2.5cm}}
 \put(6,0){\epsfig{file=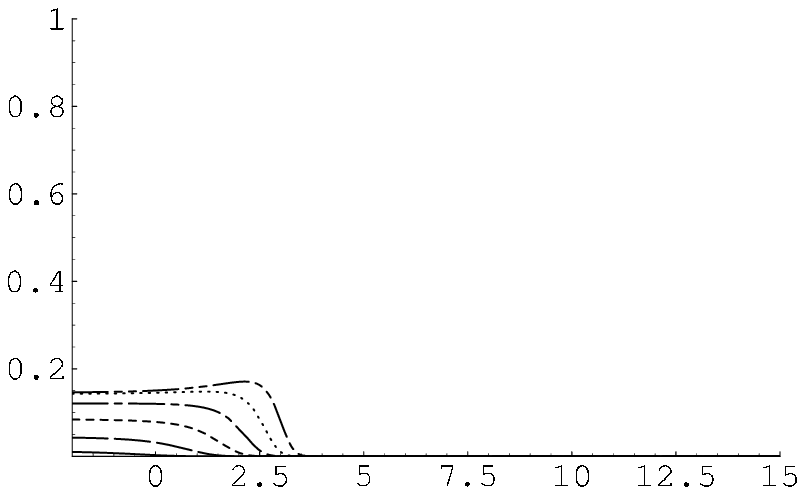,width=6cm,height=2.5cm}}
  }
 \put(-0.25,8.5){
 \marke{2,0}{cl}{P_{+,n}\,,~~\hbar\scale^2 = \frac 1 2}
 \marke{7.5,0}{cl}{P_{-,n}\,,~~\hbar\scale^2 = \frac 1 2}
 }
 \put(-0.25,4.5){
 \marke{2,0}{cl}{P_{+,n}\,,~~\hbar\scale^2 = \frac 1 8}
 \marke{7.5,0}{cl}{P_{-,n}\,,~~\hbar\scale^2 = \frac 1 8}
 }
 \put(-0.25,0.5){
 \marke{2,0}{cl}{P_{+,n}\,,~~\hbar\scale^2 = \frac 1 {32}}
 \marke{7.5,0}{cl}{P_{-,n}\,,~~\hbar\scale^2 = \frac 1 {32}}
 }
 \end{picture}

 \caption{\label{AbbPin} Expectation values $P_{+,n}$ and $P_{-,n}$
for $n=0$ (continuous line), $n=1$ (long dash), $n=2$ (short dash),
$n=3$ (long-short dash), $n=4$ (dotted line) and $n=5$ (long-short-short
dash), plotted over $\log\tinco$ for $\hbar\scale^2 = \frac 1 2,
\frac 1 8, \frac 1 {32}$.}

\end{figure}

On the level of the observables, the costratification is given by
the orthoprojectors $\Pi_\pm$ onto the subspaces $\Hi_\pm$. We
discuss their expectation values in the energy eigenstates,
$$
P_{\pm,n} := \langle \xi_n | \Pi_\pm \xi_n \rangle\,,
$$
i.e., the probability that the system prepared in the stationary state
$\xi_n$ is measured in the subspace $\Hi_\pm$. According to \eqref{Goproj},
 \beq\label{GPin}
P_{\pm,n}
  =
|\langle \xi_n | \psi_\pm \rangle|^2
\,.
 \eeq
As $\se_{2n+2}$ is odd and $\pi$-periodic, it can be expanded as
$$
\se_{2n+2}(y;q) = \sum_{k=0}\nolimits^\infty B^{2n+2}_{2k+2}(q)
\sin((2k+2)y)\,,
$$
with Fourier coefficients $B^{2n+2}_{2k+2}(q)$ satisfying certain recurrence
relations \cite[\S20.2]{AbramowitzStegun}. Due to \eqref{G-char}\,,
 \beq\label{Gscaproxichi}
\langle \xi_n | k \rangle = (-1)^{n+k}
B^{2n+2}_{2k+2}(4\tinco)\,,
 \eeq
whence \eqref{Gpsip} and \eqref{Gpsim} yield the expressions
 \beqa\label{Gxipsip}
\langle \xi_n | \psi_+ \rangle
 & = &
\frac{(-1)^n}{N}
 \sum\nolimits_{k=0}^\infty
(-1)^k\,(k+1)\,\mr e^{-\hbar\scale^2(k+1)^2/2} \,
B^{2n+2}_{2k+2}(4\tinco) ,
\\ \label{Gxipsim}
\langle \xi_n | \psi_- \rangle
 & = &
\frac{(-1)^n}{N}
 \sum\nolimits_{k=0}^\infty
(k+1)\,\mr e^{-\hbar\scale^2(k+1)^2/2} \, B^{2n+2}_{2k+2}(4\tinco) .
 \eeqa
Together with \eqref{GPin}, this procedure leads to formulas for
$P_{\pm,n}$. The functions $P_{\pm,n}$ 
depend on the parameters
$\hbar$, $\scale^2$ and $\inco$ only via the combinations
$\hbar\scale^2$ and $\tinco = \inco/(\hbar^2\scale^2)$. Figure
\rref{AbbPin} displays $P_{\pm,n}$ for $n=0,\dots,5$ as functions
of $\tinco$ for three specific values of $\hbar\scale^2$, thus
treating $\tinco$ and $\hbar\scale^2$ as independent parameters.
This is appropriate for the discussion of the dependence of
the functions $P_{\pm,n}$ 
on the coupling parameter $g$ for fixed values of
$\hbar$ and $\scale^2$. The plots have been generated by
Mathematica through numerical integration.

For $n=0$, the function $P_{+,n}$ has a dominant peak which is enclosed by less
prominent 
maxima of the other $P_{+,n}$'s and moves to higher $\tinco$ when
$\hbar\scale^2$ decreases. That is to say, for a certain value of the
coupling constant, the state $\psi_+$ which spans $\Hi_+$ seems to
coincide almost perfectly with the ground state. If the two states
coincided exactly then \eqref{Gscaproxichi} would imply that, for
a certain value of $q$, the coefficients $B^{2n+2}_{2k+2}(q)$
would be given by $(-1)^{n+k}\frac 1 N (k+1)\mr
e^{-\hbar\scale^2(k+1)^2/2}$. However, this is not true; the
latter expressions do not satisfy the recurrence relations valid
for the coefficients $B^{2n+2}_{2k+2}(q)$ for any value of $q$.

\section{Outlook}

For $\group = \SU(2)$ it remains to discuss the dynamics relative to the
costratified structure and to explore the probability flow into and out of the 
subspaces $\Hi_\pm$. More generally, it would be worthwhile carrying out this
program for $\group = \SU(n)$, $n\geq 3$. For $\group = \SU(3)$, the orbit type 
stratification of the reduced phase space consists of a 4-dimensional stratum, a
2-dimensional stratum, and three isolated points. Thereafter the approach should
be extended to arbitrary lattices. 

The notion of costratified Hilbert space implements the stratification
of the reduced classical phase space on the level of states. The significance of
the stratification for the quantum observables remains to be clarified. 
Then the physical role of this stratification can be
studied in more realistic models like the lattice QCD of
\cite{qcd1,qcd2,qcd3}.

A number of applications 
of the theory of stratified K\"ahler spaces
have already been mentioned.
Using the approach to lattice gauge theory in \cite{kan},
we intend to develop elsewhere a rigorous approach to
the quantization of certain lattice gauge theories by means of
the K\"ahler quantization scheme for 
stratified K\"ahler spaces
explained in the present paper.
We plan to apply this scheme in particular to  situations
of the kind explored in \cite{kirutwo}--\cite{kirufou} 
and to compare it with the approach to quantization
in these papers.
Constrained quantum systems occur
in molecular mechanics as well, see e.~g. \cite{taniiwai}
and the references there.
Perhaps the K\"ahler quantization scheme for 
stratified K\"ahler spaces will shed new light on these quantum systems.

\noindent 
       Johannes Huebschmann\\
       USTL, UFR de Math\'ematiques\\
       CNRS-UMR 8524\\
      59655 Villeneuve d'Ascq Cedex, France\\
      Johannes.Huebschmann@math.univ-lille1.fr
      

\label{lastpage}
\end{document}